\documentclass[11pt,reqno,oneside]{amsart}

\usepackage{graphicx}
\usepackage{subcaption}
\usepackage{amssymb}
\usepackage{epstopdf}

\usepackage[a4paper, total={6in, 9in}]{geometry}

\usepackage{amsmath,amsfonts,amsthm,mathrsfs,amssymb}
\usepackage[usenames]{color}

\usepackage{xcolor}

\usepackage{hyperref}
\usepackage{xcolor}
\hypersetup{
  colorlinks,
  linkcolor={blue!80!black},
  urlcolor={blue!80!black},
  citecolor={blue!80!black}
}
\usepackage[capitalize,noabbrev]{cleveref}

\newtheorem{thm}{Theorem}[section]

\newtheorem{lem}{Lemma}[section]

\theoremstyle{definition}
\newtheorem{defn}{Definition}[section]

\theoremstyle{remark}

\newtheorem{rem}{Remark}[section]
\numberwithin{equation}{section}

\numberwithin{equation}{section}

\newcounter{saveeqn}

\newcommand{\pare}[1]{\left(#1\right)}

\newcommand{\abs}[1]{\left\lvert #1 \right\rvert}

\newcommand{\RR}{\mathbb{R}}

\newcommand{\im}{\mathrm{i}}

\DeclareMathOperator{\curl}{curl}

\DeclareMathOperator{\cros}{\wedge}

\newcommand{\Cor}{\mathcal{K}}

\newcommand{\nor}{\nu}

\newcommand{\far}[1]{#1_{\infty}}

\newcommand{\EM}[1]{\mathbf{#1}}

\newcommand{\Ep}{\varepsilon}
\newcommand{\Mp}{\mu}
\newcommand{\Eci}{\gamma}

\title[On electromagnetic scattering from conical singularities and applications]{Visibility, invisibility and unique recovery of inverse electromagnetic problems with conical singularities }

\author{Huaian Diao}
\address{School of Mathematics, Jilin University,
Changchun, Jilin 130012, China.}
\email{diao@jlu.edu.cn}

\author{Xiaoxu Fei}
\address{School of Mathematics and Statistics, Northeast Normal University,
Changchun, Jilin 130024, China.}
\email{feixx0921@163.com}

\author{Hongyu Liu}
\address{Department of Mathematics, City University of Hong Kong, Kowloon, Hong Kong SAR, China.}
\email{hongyu.liuip@gmail.com; hongyliu@cityu.edu.hk}

\author{Ke Yang}
\address{School of Mathematics and Statistics, Northeast Normal University,
Changchun, Jilin 130024, China.}
\email{1357812761@qq.com}

\date{} 
\begin{document}
	\maketitle
	
	\begin{abstract}
		
In this paper, we study time-harmonic electromagnetic scattering in two scenarios, where the anomalous scatterer is either a pair of electromagnetic sources or an inhomogeneous medium, both with compact supports.  We are mainly concerned with the geometrical inverse scattering problem of recovering the support of the scatterer, independent of its physical contents, by a single far-field measurement. It is assumed that the support of the scatterer (locally) possesses a  conical singularity. We establish a local characterisation of the scatterer when invisibility/transparency occurs, showing that its characteristic parameters must vanish locally around the conical point. Using this characterisation, we establish several local and global uniqueness results for the aforementioned inverse scattering problems, showing that visibility must imply unique recovery. In the process, we also establish the local vanishing property of the electromagnetic transmission eigenfunctions around a conical point under the H\"older regularity or a regularity condition in terms of Herglotz approximation.

		\medskip

		\noindent{\bf Keywords:}~~electromagnetic waves, geometrical inverse scattering, conical singularity, invisibility and transparency, locally vanishing, unique recovery, single far-field measurement, transmission eigenfunctions. 
		
		\noindent{\bf 2010 Mathematics Subject Classification:}~~78A45, 35Q61, 35P25 (primary); 78A46, 35P25, 35R30 (secondary).
		
	\end{abstract}

	\section{Introduction}
	
	In this paper, we study time-harmonic electromagnetic scattering in two scenarios, where the anomalous scatterer is either a pair of electromagnetic sources or an inhomogeneous medium, both with compact supports. We first introduce the forward scattering problems in the two scenarios. 

Let $\Omega$ be a bounded Lipschitz domain in $\mathbb R^3$ with a connected complement $\mathbb{R}^3\backslash\overline{\Omega}$, which signifies the support of an inhomogeneous scatterer. Let $\mathbf{J}_j(\mathbf{x})$, $\mathbf{x}\in\mathbb{R}^3$ and $j=1,2$, be $\mathbb{C}^3$-valued functions such that $\mathrm{supp}(\mathbf{J}_j)\subset\Omega$. It is assumed that $\mathbf{J}_1|_{\Omega}\in L^2(\Omega; \mathbb{C}^3)$ and $\mathbf{J}_2|_{\Omega}\in L^2(\Omega; \mathbb{C}^3)$, which signify the intensities of active electric and magnetic sources, respectively. Let $(\mathbf{E},\mathbf{H})\in H_{\mathrm{loc}}(\mathrm{curl}, \mathbb{R}^3)\times H_{\mathrm{loc}}(\mathrm{curl}, \mathbb{R}^3)$ denote the electric and magnetic fields respectively. 
Consider the following  electromagnetic  scattering problem 
	\begin{equation}\label{eq:Maxwell1}
		\begin{cases}
			&\nabla\wedge \EM{E}(\mathbf x)-\im\omega \Mp_0 \EM{H}(\mathbf x)=\EM{J}_1(\mathbf x),\quad \mathbf x\in\RR^3,\\
			& \nabla\wedge \EM{H}(\mathbf x)+\im\omega\Ep_0 \EM{E}(\mathbf x)=\EM{J}_2(\mathbf x),\quad\, \mathbf x\in\RR^3,\\
			& \lim_{|\mathbf x|\to\infty} |\mathbf x|\pare{{\mu}_0^{1/2}  \EM{H}\times\frac{\mathbf x}{|\mathbf x|}-{\varepsilon}_0^{1/2}\EM{E}}=0,	
		\end{cases}
	\end{equation}
	where $\omega\in\mathbb{R}_+$ signifies the frequency of the wave, and $\varepsilon_0$ and $ \mu_0\in\mathbb{R}_+$, respectively, denote the electric permittivity and magnetic permeability of a
	uniformly homogeneous space. The last limit in \eqref{eq:Maxwell1} is known as the Silver-M\"uller radiation condition which holds uniformly in all directions
	$\hat {\mathbf x}:=\mathbf x/|\mathbf x|\in\mathbb{S}^2$, $\mathbf x\in\mathbb{R}^3\setminus\{\mathbf 0\}$, and characterizes the outgoing nature of the electromagnetic waves. It is emphasized that we consider the possible presence of both electric and magnetic sources, though only the electric source might be the physically meaningful one. The well-posedness of the Maxwell system \eqref{eq:Maxwell1} can be conveniently found in \cite{LRX16,Ned01}. We know that as $|{\mathbf x}|\rightarrow+\infty$ it holds that    
	\begin{equation}\label{eq:radiation}
	\pare{\EM{E},\EM{H}}(\mathbf x)=\frac{e^{\im k |\mathbf x|}}{|\mathbf x|} \pare{\far{\EM{E}},\far{\EM{H}}}(\hat{\mathbf x})+\mathcal{O}\left(\frac{1}{|\mathbf x|^2}\right), 
	\end{equation}
	where $k:=\omega\sqrt{\varepsilon_0\mu_0}$ is known as the wavenumber, and $\mathbf{E}_\infty(\hat {\mathbf x})$ and $\mathbf{H}_\infty(\hat {\mathbf x})$ are referred to the electric and the magnetic far-field patterns respectively. It is known that $\mathbf{E}_\infty$ and $\mathbf{H}_\infty$ are analytic functions on the unit sphere $\mathbb{S}^2$ with the following one-to-one correspondence
	\begin{equation}\label{eq:equiv1}
	\mathbf{H}_\infty(\hat {\mathbf x})=\hat {\mathbf x}\wedge \mathbf{E}_\infty(\hat {\mathbf x})
	\quad\text{and}\quad
	\mathbf{E}_\infty(\hat {\mathbf x})=-\hat {\mathbf x}\wedge \mathbf{H}_\infty(\hat {\mathbf x}),\quad\ \ \forall \hat {\mathbf x}\in\mathbb{S}^2. 
	\end{equation}

Next we introduce the scattering due to the interaction of a (passive) inhomogeneous medium scatterer and an (actively sent) incident wave. Suppose that $\Omega$ supports an inhomogeneous medium whose material parameters are characterised by 
 the electric permittivity $\varepsilon\in L^\infty(\Omega; \mathbb{R}_+)$, magnetic permeability $\mu\in L^\infty(\Omega; \mathbb{R}_+)$ and electric conductivity $\sigma\in L^\infty(\Omega; \mathbb{R}_+^0)$. Throughout the rest of the paper,  we set
	$$
	\varepsilon(\mathbf x)=\begin{cases}
		\varepsilon_0, \ \mathbf x\in \mathbb{R}^3\setminus\overline{\Omega}, \\
		\varepsilon,\ \ \mathbf x\in \Omega,
	\end{cases},\ \mu(\mathbf x)=\begin{cases}
		\mu_0, \  \mathbf x\in \mathbb{R}^3\setminus\overline{\Omega}, \\
		\mu,\ \ \mathbf x\in \Omega,
	\end{cases},\ \sigma(\mathbf x)=\begin{cases}
		0, \ \mathbf x\in \mathbb{R}^3\setminus\overline{\Omega}, \\
		\sigma,\  \mathbf x\in \Omega.
	\end{cases} 
	$$
The incident wave field $(\mathbf{E}^i, \mathbf{H}^i)$  is a pair of entire solutions to 
\begin{equation}\label{eq:Maxwellh}
		\nabla \wedge  \EM{E}^i-\im\omega \Mp_0 \EM{H}^i=\mathbf 0,\quad \nabla\wedge \EM{H}^i+\im\omega\Ep_0 \EM{E}^i=\mathbf 0\quad\mbox{in}\ \ \mathbb{R}^3,
	\end{equation}	
which interacts with the scattering medium described above. The resulting scattered electromagnetic  wave field is denoted by $(\mathbf{E}, \mathbf{H})$. Hence the total electromagnetic  wave field  $(\mathbf{E}^t, \mathbf{H}^t)$ is the superposition of the incident and scattered waves, i.e.,  
\[
	(\mathbf{E}^t, \mathbf{H}^t)=(\mathbf{E}^i, \mathbf{H}^i)+(\mathbf{E}, \mathbf{H})\quad\mbox{in}\ \ \mathbb{R}^3. 
	\]
The aforementioned electromagnetic  scattering  can be modelled by the following  Maxwell system 
	\begin{equation}\label{eq:Maxwell2}
		\begin{cases}
			& \nabla\wedge\mathbf{E}^t({\mathbf x})-\mathrm{i}\omega\mu({\mathbf x}) \mathbf{H}^t({\mathbf x})=0,\hspace*{1.8cm} {\mathbf x}\in\mathbb{R}^3,\medskip\\
			&\nabla\wedge\mathbf{H}^t({\mathbf x})+\mathrm{i}\omega\gamma({\mathbf x})\mathbf{E}^t({\mathbf x})=0,\hspace*{1.8cm}{\mathbf x}\in\mathbb{R}^3,\medskip\\
			& \lim_{|\mathbf x|\to\infty} |\mathbf x|\pare{{\mu}_0^{1/2}\EM{H}\times\frac{\mathbf x}{|\mathbf x|}-{\varepsilon}_0^{1/2}\EM{E}}=0. 
		\end{cases}
	\end{equation}
where $\gamma(\mathbf x)=\varepsilon(\mathbf x)+\mathrm{i}\sigma(\mathbf x)/\omega$. The well-posedness of \eqref{eq:Maxwell2} can be found in \cite{LRX16,Ned01}, which guarantees the unique existence of a pair of solutions
	$(\mathbf{E},\mathbf{H})\in H_{\mathrm{loc}}(\mathrm{curl}, \mathbb{R}^3)\times H_{\mathrm{loc}}(\mathrm{curl}, \mathbb{R}^3)$ to  \eqref{eq:Maxwell2}. Furthermore, the scattered electromagnetic  waves $\pare{\EM{E},\EM{H}}$ have the following asymptotic expansion: 
	\begin{equation}\label{eq:farfield}
		\pare{\EM{E},\EM{H}}(\mathbf x)=\frac{e^{\im k |\mathbf x|}}{|\mathbf x|} \pare{\far{\EM{E}},\far{\EM{H}}}(\hat{\mathbf x})+\mathcal{O}\left(\frac{1}{|\mathbf x|^2}\right). 
	\end{equation}

Henceforth, we write $(\Omega; \mathbf{J}_1, \mathbf{J}_2)$ and $(\Omega; \varepsilon, \mu, \sigma)$ to denote the source and medium scatterers respectively introduced above. Here, $\Omega$ signifies the support of the scatterer which contains its shape and location information, whereas $\mathbf{J}_1, \mathbf{J}_2$ or $\varepsilon, \mu, \sigma$ are its physical content, and hereafter are referred to as the characteristic parameters of the scatterer. In the case that $\Omega$ is a medium scatterer, we also include the incident field $(\mathbf{E}^i, \mathbf{H}^i)$ as characteristic parameter since its interaction with the medium parameters generates the source that produces the radiating scattered waves. In this paper, one of the major concerns is the following geometrical inverse scattering problem: 
\begin{equation}\label{eq:ip2}
		\mathbf{E}_\infty(\hat {\mathbf x}),\ \hat {\mathbf x}\in\mathbb{S}^2 \mbox{ or }  \mathbf{H}_\infty(\hat {\mathbf x}),\ \hat {\mathbf x}\in\mathbb{S}^2  \longmapsto \Omega\ \ \mbox{independent of its physical content}, 
	\end{equation}
where $\mathbf{E}_\infty$ (or, equivalently $\mathbf{H}_\infty$ by virtue of \eqref{eq:equiv1}) is either from \eqref{eq:radiation} for the source scattering or \eqref{eq:farfield} for the medium scattering. It is straightforwardly verified that the inverse problem \eqref{eq:ip2} is nonlinear, though the forward scattering problem is linear. Throughout our study, it is assumed for \eqref{eq:ip2} that $\omega\in\mathbb{R}_+$ is fixed and in case $\Omega$ is a medium scatterer, the far-field pattern in \eqref{eq:ip2} is collected corresponding to a single incident wave field $(\mathbf{E}^i, \mathbf{H}^i)$. In such a case, $\mathbf{E}_\infty$ is referred to as a single far-field measurement. Clearly, in order to determine $\Omega$, it is sufficient to recover $\partial\Omega$. It can be seen that the inverse problem \eqref{eq:ip2} is formally determined with a single far-field measurement since both $\mathbb{S}^2$ and $\partial\Omega$ are two-dimensional manifolds. The geometrical inverse problem \eqref{eq:ip2} is a well-known longstanding one in the inverse scattering theory \cite{CoK13,Liu22,LZ06}, with a colourful history and yet still largely open. It lays the theoretical foundation for many wave imaging technologies including radar, medical imaging and non-destructive testing where one is more interested in extracting the geometrical information of the anomalies by limited measurement data.

In respect to \eqref{eq:ip2}, a closely related problem is the occurrence of invisibility/transparency, namely $\mathbf{E}_\infty=\mathbf{H}_\infty\equiv \mathbf{0}$. In such case, $(\Omega; \mathbf{J}_1, \mathbf{J}_2)$ is said to be a non-radiating/raditioneless  source and $(\Omega; \varepsilon, \mu, \sigma)$ is said to be a transparent/invisible scatterer. Recently, geometrical characterisations of non-radiating sources and transparent/invisible mediums have received considerable interest in the literature; see \cite{Bsource,BL2016,BL2017,BPS,CV,CaX,DDL1,DCL,ElH,PSV,SS,VX} for related studies in acoustic scattering,  \cite{BLY,DLS} in elastic scattering and \cite{BLX2020,2021,LX}  in electromagnetic scattering. Roughly speaking, if the support of the scatterer possesses a certain geometrical singularity on its boundary $\partial\Omega$, say e.g. a corner, then the characteristic parameters of an invisible scatter must be vanishing (locally) around the geometrical singular point. As a direct consequence, if the characteristic parameters of a scatterer is a-priori known to be non-vanishing around a geometrically singular point, then it must radiate a nontrivial scattering pattern, namely it must be visible with respect to far-field measurement. It is emphasized that this point has been essentially implied in all of the aforementioned studies on characterising radiating/non-radiating scatterers, though it may appear in different phrasings. It is also noted that in \cite{BL2018}, a smooth boundary point with a sufficiently high curvature is shown to possess a similar characterisation as above in the context of acoustic scattering. In the current article, we make a novel contribution along the line by establishing the local vanishing property for the electromagnetic scattering in the two scenarios introduced earlier when the scatterer possesses a conical singularity. It is noted that in the context of electromagnetic scattering, only polyhedral singularities have been considered due to the highly complicated physical and technical nature. The main results on this aspect are contained in Theorems~
\ref{thm:thm2.2} and \ref{thm:main5}, respectively, for the source and medium scattering.  

\begin{figure}
	\centering
\begin{subfigure}{0.4\textwidth}\includegraphics[width=0.8\linewidth]{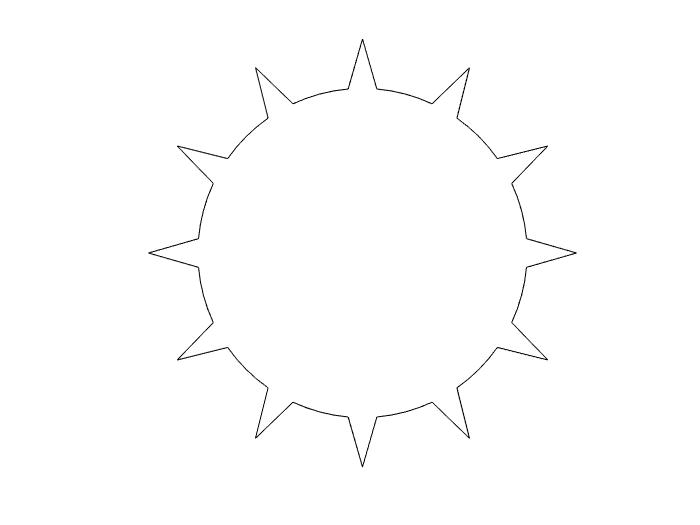}\end{subfigure}\hspace*{-2.1cm}\begin{subfigure}{0.4\textwidth}\includegraphics[width=1\linewidth]{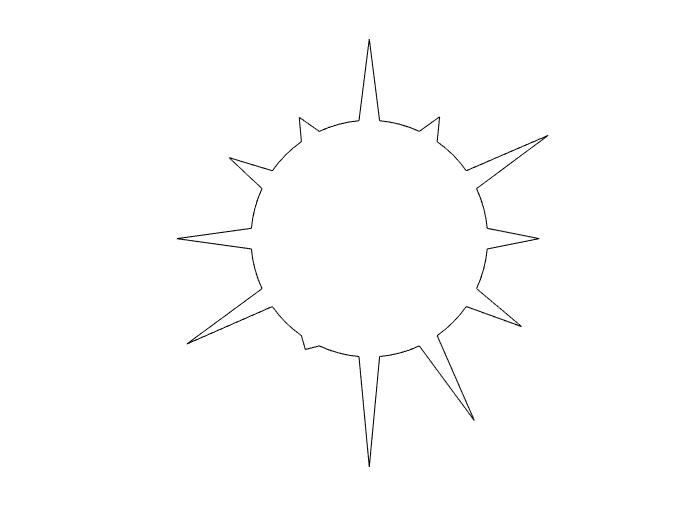}\end{subfigure}\hspace*{-.5cm}\begin{subfigure}{0.4\textwidth}\vspace*{-.2cm}\includegraphics[width=0.8\linewidth]{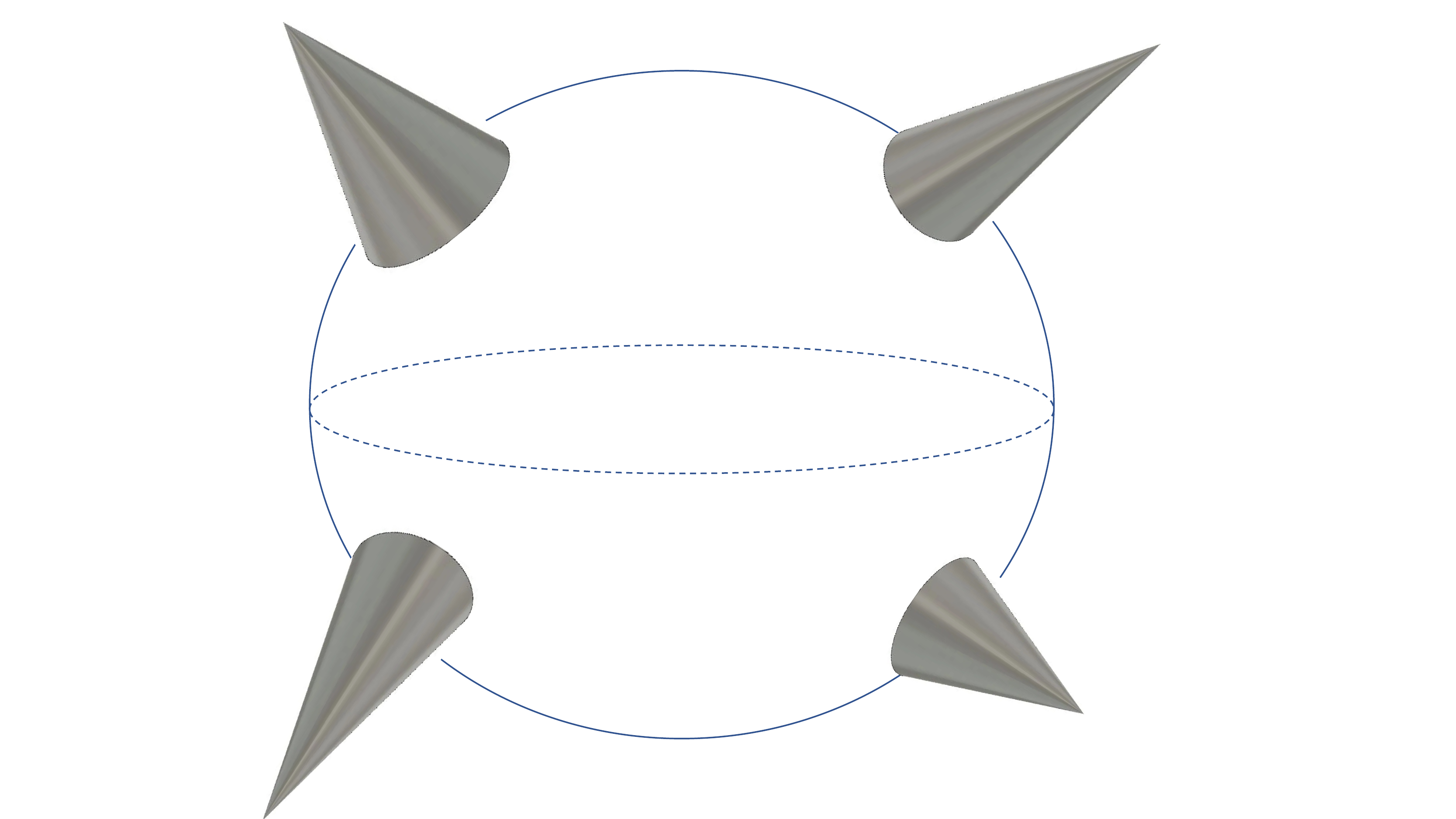}\end{subfigure}
	\caption{Schematic illustration of coronal-shape scatterers. Rigorous definition is provided in Definition~\ref{def:finite c}. The first two are the slice plottings of two coronal-shape scatterers with many conical singularities and the third one is a 3D plotting with 4 conical singularities on its body.}
	\label{fig:coordinate1}
\end{figure}

If visibility is guaranteed, namely the scatterer does generate scattering information to the far-field observer, the next issue of primary importance to \eqref{eq:ip2} is the unique identifiability. That is, if there are two scatterers $\Omega$ and $\Omega'$, with possibly different and not a-priori known physical contents, which generate the same far-field measurement if and only if $\Omega=\Omega'$. By using the geometrical characterisation discussed above for non-radiating/transparent scatterers, we establish several novel local and global unique recovery results for \eqref{eq:ip2}, showing that visibility is equivalent to unique recovery. It is clear that for the inverse problem \eqref{eq:ip2}, visibility, i.e. $\mathbf{E}_\infty$ is not identically zero, does not necessarily imply unique identifiability. In our study, we can achieve such an equivalence relation due to the fact that our analysis is localised around the conical point. It is interesting to note that our global recovery results contain a special case that the scatterer is of coronal shape (cf. Fig.~\ref{fig:coordinate1} for a schematic illustration), which may be of practical interest to the medical imaging. Finally, we would like to mention in passing some related results on uniqueness for geometrical inverse electromagnetic problems by a single far-field measurement \cite{BLX2020,DLZZ21,LRX16,LYZ}. 

Finally, we also achieve a geometrical characterisation of electromagnetic transmission eigenfunctions, showing that they must vanish (locally) around a conical point. Transmission eigenvalue problems arise from non-scattering/invisibility but go beyond, especially when the regularity of transmission eigenfunctions is weakened; see Section~\ref{sec:4}  for more related background discussion. Recently, the spectral geometry of transmission eigenfunctions has also received considerable attention in the literature; see \cite{Bsource,BLY,BL2017b,BL2018,BLLW,CaX,DCL,DLS} in different physical context and especially \cite{BLX2020,2021} in the context of electromagnetic scattering for local structures and \cite{CDHLW,DJLZ,DLX1} for global structures. We establish a local vanishing property of the electromagnetic transmission eigenfunctions around a conical point under the H\"older regularity or a regularity condition in terms of Herglotz approximation, which add a novel contribution to the spectral theory of transmission eigenfunctions. 

According to our discussion above, the visibility, invisibility and unique recovery of the inverse electromagnetic problem \eqref{eq:ip2} as well as spectral geometry of transmission eigenfunctions are separate but intriguingly connected topics. We present all those geometrical results as discussed above to corroborate the interesting connections among them. Finally, we would like to briefly discuss the mathematical strategy in establishing those geometrical results. We shall make essential use of tools from microlocal analysis to carefully analyse the singularity behaviour of the solution to the Maxwell system induced by the geometrical singularity of the shape of the underlying scatterer. This shares a similar spirit to \cite{BLX2020} which deals with a polyhedral corner. Nevertheless, we achieve several new technical developments in order cope with the different geometrical setup as well as several other issues, especially to significantly weaken the regularity assumptions needed in \cite{BLX2020}, and make the study more physically relevant.

The rest of the paper is organised as follows. In Section \ref{sec:2}, we consider the geometrical characterisation of non-radiating sources as well as the geometrical inverse problem \eqref{eq:ip2} in determining the support of a source scatterer. In Section \ref{sec:3}, we consider the geometrical characterisation of transparent/invisible medium scatterers as well as the geometrical inverse problem \eqref{eq:ip2} in determining the support of a medium scatterer. Section \ref{sec:4} is devoted to the geometrical characterisation of electromagnetic transmission eigenfunctions.

\section{Non-radiating sources and inverse source scattering}\label{sec:2}

In this section, we establish the vanishing property of non-radiating electromagnetic sources around a conical corner, and then use it to derive unique shape determination  for the inverse electromagnetic source problem.  We first introduce the geometric setup of our study.
	
Let $\mathbf{a}\in\mathbb{S}^2:=\{\mathbf{x}\in\mathbb{R}^3~|~|\mathbf{x}|=1\}$, $\mathbf{x}_0\in\mathbb{R}^3$ and $\theta_0\in (0, \pi/2)$ be fixed. Define
\begin{equation}\label{eq:ccone}
\Cor_{\mathbf x_0,\theta_0}:=\left\lbrace \mathbf{x}=\mathbf{x}_0+r\mathbf{\tilde{x}}~|~\langle \mathbf{\tilde{x}}, \mathbf{a}\rangle\in (0, \theta_0), \forall \mathbf{\tilde{x}}\in\mathbb{S}^2,\ \ \forall r\in\mathbb{R}_+\right\rbrace, 
\end{equation}
$ \Cor $ is a convex cone with an opening angle $ 2\theta_{0} $ less than $ \pi $, where  $ \mathbf x_{0} $ is the apex of the cone and $\mathbf a$ is the axis of $ \Cor $. Given a constant $r_0\in\mathbb{R}_+$, we define
	\begin{equation}\label{eq:coner0}
	\Cor_{r_0}=	\Cor^{r_0}_{\mathbf x_0}=\mathcal \Cor \cap B_{r_0}(\mathbf x_0),
	\end{equation}
	where  $ B_{r_0}(\mathbf x_0):=\{\mathbf x \in\mathbb R^3~|~ |\mathbf x-\mathbf x_0 | < r_0  \}$. Without loss of generality, throughout this paper, we let $\mathbf x_0$ be the origin and the axis  $ \mathbf a =\mathbf e_3$ with $\mathbf e_3=(0,0,1)^\top $.

\subsection{Geometrical characterisation of non-radiating sources}

In order to prove the geometrical characterisation of a radiationless electromagnetic source near a  conical corner, we need the following  lemmas. 
	\begin{lem}\label{lem:Integral}\cite[Lemma 2.1]{BLX2020}
		Let $\Omega$ be a bounded Lipschitz domain in $\RR^3$ and $\EM{J}_j\in L^2(\Omega; \mathbb{C}^3)$.
		Suppose that $\pare{\EM{E},\EM{H}}\in H(\mathrm{curl},\Omega)\times H(\mathrm{curl},\Omega)$ is a solution to the Maxwell system 
		\begin{equation}\label{eq:MaxwellO}
		\nabla\wedge  \EM{E}-\im\omega \Mp_0 \EM{H}=\EM{J}_1,\quad \nabla\wedge \EM{H}+\im\omega\Ep_0 \EM{E}=\EM{J}_2\qquad \mbox{in $\Omega$}.
		\end{equation}
		Then one has
		\begin{equation}\label{eq:IntId1}
		\int_{\Omega} \EM{J}_1\cdot\EM{W} \rm d \mathbf x 
		+\int_{\Omega} \EM{J}_2\cdot\EM{V} \rm d \mathbf x 
		=\int_{\partial \Omega} \EM{W}\cdot\pare{\nor\cros\EM{E}} \rm d \sigma
		+\int_{\partial \Omega} \EM{V}\cdot\pare{\nor\cros\EM{H}} \rm d \sigma ,
		\end{equation}
		and
		\begin{equation}\label{eq:IntId2}
		\Ep_0\int_{\Omega} \EM{J}_1\cdot\EM{V} \rm d \mathbf x 
		-\Mp_0\int_{\Omega} \EM{J}_2\cdot\EM{W} \rm d \mathbf x 
		=\Ep_0\int_{\partial \Omega} \EM{V}\cdot\pare{\nor\cros\EM{E}}\rm d \sigma
		-\Mp_0\int_{\partial \Omega} \EM{W}\cdot\pare{\nor\cros\EM{H}} \rm d \sigma ,
		\end{equation}
		for any $\pare{\EM{V},\EM{W}}\in H(\mathrm{curl},\Omega)\times H(\mathrm{curl},\Omega)$ satisfying 
		\begin{equation}\label{eq:MaxwellTest}
		\nabla\wedge  \EM{V}-\im\omega \Mp_0 \EM{W}=0,\quad \nabla\wedge \EM{W}+\im\omega\Ep_0 \EM{V}=0\quad \mbox{ in } \quad \Omega.
		\end{equation}
	\end{lem}
	
	From the proof of Theorem 1.1 in \cite{BLX2020}, we can summarize the following lemma:
	
	\begin{lem}
		For any given vectors $ \mathbf d, \mathbf d^{\perp}\in \mathbb S^2 $ such that $ \mathbf d^{\perp} \perp \mathbf d $, denote the complex vectors
		\begin{equation}\label{eq:r}
		\boldsymbol{\rho}= \tau \mathbf d + \im \sqrt{\tau^2+k^2} \mathbf d^\perp,\ \mathbf p=\mathbf d^\perp- \im \sqrt{1+k^2/\tau^2} \mathbf d,
		\end{equation}
		where $\tau\in \mathbb R_+$ and $ k=\omega \sqrt{\Ep_0\Mp_0} $ with $\omega, \Ep_0, \Mp_0 \in \mathbb R_+$, 
		then one has 
		\begin{equation}\label{eq:p norm}
		\boldsymbol{\rho} \cdot \mathbf p=0,\quad \boldsymbol{\rho} \cros \mathbf p=-k^2\left( \mathbf d \cros \mathbf d^\perp\right) /\tau.
		\end{equation}
		Let 
		\begin{equation}\label{eq:vw}
		\EM{V}(\mathbf x)=\mathbf pe^{\boldsymbol{\rho}\cdot \mathbf  x}\quad\text{and}\quad 
		\EM{W}(\mathbf x)=\frac{1}{\im \omega\mu_0} \boldsymbol{\rho} \cros \mathbf p e^{\boldsymbol{\rho}\cdot \mathbf x},
		\end{equation}
		then $\pare{\EM{V},\EM{W}}$ is a pair of solutions to the Maxwell system \eqref{eq:MaxwellTest}. Futhermore, if we choose 
		\begin{equation}\label{eq:vw2}
			\EM{V}(\mathbf x)=-\frac{1}{\im \omega\Ep_0} \boldsymbol{\rho}\cros \mathbf p e^{\boldsymbol{\rho}\cdot \mathbf x}\quad\text{and}\quad 
			\EM{W}(\mathbf x)=\mathbf pe^{\boldsymbol{\rho}\cdot \mathbf x},
		\end{equation}
		$\pare{\EM{V},\EM{W}}$ is also a pair of solutions to the Maxwell system \eqref{eq:MaxwellTest}.
	\end{lem}

	

	\begin{lem}\cite[Lemma 4.4]{DLW21}\label{lem:infty}
		Let $ \alpha>0 $ and $ 0<\delta<e $ be two given parameters, $ \eta \in \mathbb C  $ and $ \Re (\eta)>0 $, we have
		\begin{equation}\label{eq:epis}
		\int_{0}^{\delta}r^{\alpha}e^{-\eta t} {\mathrm d} t=\Gamma(\alpha +1)/\eta^{\alpha +1}-\int_{\delta}^{\infty}r^{\alpha}e^{-\eta r}{\mathrm d}r.
		\end{equation}  
		When $ \Re \eta\ge  \frac{2\alpha}{e} $, it yields that
		\begin{equation}\notag
		\left| \int_{\delta}^{\infty}r^{\alpha}e^{-\eta r}{\mathrm d}r  \right| \le \frac{2}{\Re \eta}e^{-\frac{\delta}{2}\Re \eta}.
		\end{equation}
	\end{lem}
		
		
		
	 
In the following lemma, we shall establish a key asymptotic analysis of the integral \eqref{eq:NonvanishIntegral} with respect to the parameter $ \tau  $ goes to infinity, which shall play an important role in proving Theorem~\ref{thm:main0}.


	\begin{lem}\label{lem:Nonvanish}
		Let $ \Cor $ be defined by \eqref{eq:ccone} with the apex at the origin, where the opening angle of $ \Cor $ is $ 2\theta_{0} \in (0,{\pi}) $. 
		Let  $\Cor_{r_{0}}=\Cor\cap B_{r_{0}}$,  where  $ B_{r_0}:=\{\mathbf x \in\mathbb R^3~|~ |\mathbf x | < r_0  \}$. Then there exist a positive constant $ \delta$ and vectors $ \mathbf d, \mathbf d^{\perp}\in \mathbb S^2 $, $ \mathbf d \perp \mathbf d^{\perp} $ satisfying:
		\begin{equation}\label{eq:innerdot}
		\mathbf d \cdot \hat{\mathbf x} \le -\delta,\quad  \forall \hat{\mathbf x}\in \Cor \cap \mathbb S^2. 
		\end{equation}
		Suppose that  $\boldsymbol{\rho}$ is given by \eqref{eq:r}, where $\mathbf d$ fulfills \eqref{eq:innerdot}. Then
		\begin{equation}\label{eq:NonvanishIntegral}
		\abs{\int_{\Cor_{r_{0}}} e^{\boldsymbol{\rho}\cdot \mathbf x} {\mathrm d} \mathbf x} 
		\ge C_{\Cor} \frac{\tau^{-3}}{(1+\frac{k^2}{\tau^2})^{3/2}}+\mathcal{O}(\tau^{-1}e^{-\frac{1}{2}r_{0}\tau \delta }),
		\end{equation}
		holds for $\tau$ sufficiently large,  
		where  $ C_{\Cor}=\sqrt{2}\pi(1-\cos\theta_{0}) $ is a positive constant. 
	\end{lem}
	\begin{proof}
		Using polar coordinates transformation and Lemma~\ref{lem:infty}, we can deduce that
		\begin{equation}\notag
		\begin{split}
		\int_{\Cor_{r_{0}}} e^{\boldsymbol{\rho}\cdot \mathbf x} {\mathrm d} \mathbf x
		=&\int_{0}^{2\pi} {\mathrm d} \varphi \int_{0}^{\theta_{0}}(\frac{\Gamma(3)}{\tau^3(\mathbf d\cdot \hat{\mathbf x}+\im \sqrt{1+(\frac{k}{\tau })^2}\mathbf d^{\perp }\cdot \hat{\mathbf x})^3}+I_{r}) \sin\theta {\mathrm d} \theta:=I_{1}+I_{2},
		\end{split}		
		\end{equation}
		where
		\begin{equation}\notag
		\begin{split}
		&I_{1}=\int_{0}^{2\pi}\int_{0}^{\theta_{0}}(\frac{\Gamma(3)}{\tau^3(\mathbf d\cdot \hat{\mathbf x}+\im \sqrt{1+(\frac{k}{\tau })^2}\mathbf d^{\perp }\cdot \hat{\mathbf x})^3}+I_{r}) \sin\theta {\mathrm d} \varphi {\mathrm d} \theta,\\ &I_{2}=\int_{0}^{2\pi}\int_{0}^{\theta_{0}}I_{r}{\mathrm d} \varphi {\mathrm d} \theta,\quad 
		I_{r}=\int_{r_{0}}^{+\infty}r^{2}e^{-\tau r(\mathbf d\cdot \hat{\mathbf x}+\im \sqrt{1+(\frac{k}{\tau })^2}\mathbf d^{\perp }\cdot \hat{\mathbf x})}{\mathrm d}r.
		\end{split}
		\end{equation} 
	Using the integral mean value theorem, it  yields that
		\begin{equation}\notag
		\begin{split}
		I_{1}=&\frac{2}{\tau^3}\int_{0}^{2\pi}\frac{1}{(\mathbf d\cdot \hat{\mathbf x}(\varphi,\theta_{\xi})+\im \sqrt{1+(\frac{k}{\tau })^2}\mathbf d^{\perp }\cdot \hat{\mathbf x}(\varphi,\theta_{\xi}))^3} {\mathrm d} \varphi \int_{0}^{\theta_{0}}\sin\theta {\mathrm d} \theta\\
		=&\frac{2\cdot 2\pi (1-\cos \theta_{0})}{\tau^3}\dfrac{1}{(\mathbf d\cdot \hat{\mathbf x}(\varphi_{\xi},\theta_{\xi})+\im \sqrt{1+(\frac{k}{\tau })^2}\mathbf d^{\perp }\cdot \hat{\mathbf x}(\varphi_{\xi},\theta_{\xi}))^3},\\
		\end{split}
		\end{equation}
		where $ \theta_{0} \in (0,\frac{\pi}{2})$. Due to Lemma~\ref{lem:infty}, we can obtain 
		\begin{equation}\notag
		\left|I_{r} \right| \le \frac{2}{\tau \mathbf d\cdot \hat{\mathbf x}}e^{\frac{r_{0}}{2}\tau \mathbf d\cdot \hat{\mathbf x}}, 
		\end{equation}
which can be used to deduce that
		\begin{equation}\notag
		\left| I_{2} \right| \le \int_{0}^{2\pi }\int_{0}^{\theta_{0}}\left|I_{r} \right| {\mathrm d} \varphi {\mathrm d} \theta \le  \int_{0}^{2\pi }\int_{0}^{\theta_{0}} \frac{2}{\tau \mathbf d\cdot \hat{\mathbf x}}e^{\frac{r_{0}}{2}\tau \mathbf d\cdot \hat{\mathbf x}} {\mathrm d} \varphi {\mathrm d} \theta. 
		\end{equation}
By virtue of \eqref{eq:innerdot}, one has 
		\begin{equation}\notag
		\begin{split}
		\left| \int_{\Cor_{r_{0}}} e^{\boldsymbol{\rho}\cdot \mathbf x} {\mathrm d} \mathbf x \right| 
		\ge &\frac{4\pi(1-\cos\theta_{0})}{\tau^3}\cdot \frac{1}{2^{3/2}(1+\frac{k^2}{\tau^2})^{3/2}}+\frac{4\pi \theta_{0} }{\tau \delta }e^{-\frac{1}{2}r_{0}\tau \delta }\\
		=&C_{\Cor} \frac{\tau^{-3}}{(1+\frac{k^2}{\tau^2})^{3/2}}+\frac{4\pi \theta_{0} }{\tau \delta }e^{-\frac{1}{2}r_{0}\tau \delta },
		\end{split}
		\end{equation}
		which readily implies \eqref{eq:NonvanishIntegral} for $\tau$ sufficiently large. 
		
		The proof is complete. 
	\end{proof}

	
	From the proof of Theorem 2.8 in \cite{2021}, we have the following lemma.  
	\begin{lem}\label{lem:vani yang}
		Let $\mathbf d=(0,0,-1)^\top$ and $\mathbf d^{\perp}=(\cos\varphi ,\sin \varphi , 0)^\top$,  where $\varphi\in (0,2\pi]$. Suppose that $\mathbf p$ is defined in \eqref{eq:r}, where $\tau$ is a positive parameter of $\mathbf p$ in \eqref{eq:r}. For any complex vector $\mathbf a \in \mathbb C^3$, if 
		$$
		\lim_{\tau \rightarrow +\infty } \mathbf a \cdot \mathbf p =0, 
		$$
		then $\mathbf a=\mathbf 0$. 
	\end{lem}
	
	In order to prove the vanishing property of a non-radiating electric  and magnetic sources near a conical corner, we first need the following theorem, which shall also be used to prove geometrical characterization of a transparent/invisible medium and the vanishing of electromagnetic transmission eigenfunctions around a conical corner in what follows.

	\begin{thm}\label{thm:main0}
Suppose that $\mathbf{J}_1\in C^\alpha (\overline{\Cor^{r_0}_{\mathbf x_0}})^3$ and $\mathbf{J}_2\in C^\alpha  (\overline{ \Cor^{r_0}_{\mathbf x_0}})^3$, where  $\Cor^{r_0}_{\mathbf x_0}$	is defined by \eqref {eq:coner0} and $\alpha \in (0,1)$.   Consider the following time-harmonic electromagnetic system: 
		\begin{equation}\label{eq:MaxwellCor0}
		\left\{
		\begin{split}
		& \nabla\wedge  \EM{E}-\im\omega \Mp_0 \EM{H}=\EM{J}_{1} \quad \mbox{in $\Cor^{r_0}_{\mathbf x_0}$},\\
		& \, \nabla\wedge \EM{H}+\im\omega\Ep_0 \EM{E}=\EM{J}_{2} \quad \mbox{in $\Cor^{r_0}_{\mathbf x_0}$},\\
		&\, \nor\cros \EM{E}=\nor\cros\EM{H}=\mathbf 0 \hspace*{0.8cm}\mbox{on $\partial \Cor^{r_0}_{\mathbf x_0}\setminus \partial B_{r_0}(\mathbf x_0)$}, 
		\end{split}\right.
		\end{equation}
		where  $ (\mathbf E, \mathbf H)\in H({\rm curl},\Cor_{\mathbf x_{0}}^{r_0}) \times H({\rm curl},\Cor_{\mathbf x_{0}}^{r_0}) $, and $\nu\in\mathbb{S}^2$ is the exterior unit normal vector to $\partial \Cor^{r_0}_{\mathbf x_0}\setminus \partial B_{r_0}(\mathbf x_0)$. Then it holds that
				\begin{equation}\label{eq:cond2 thm21}
	\mathbf{J}_1({\mathbf x}_0)=\mathbf{J}_2({\mathbf x}_0)=\mathbf  0. 
	\end{equation}	
	\end{thm}

	\begin{proof}
		Since the operator $\nabla \wedge$ is invariant under rigid motion, without loss of generality, we assume that $\mathbf x_0=\mathbf 0$.  By virtue of Lemma~\ref{lem:Integral} and the boundary condition in \eqref{eq:MaxwellCor0}, for any $\pare{\EM{V},\EM{W}}$ satisfying \eqref{eq:MaxwellTest}, there holds
		\begin{equation}\label{eq:IntIdProof}
		\int_{\Cor_{r_{0}}} \EM{J}_1\cdot\EM{W}\mathrm d \mathbf x
		+\int_{\Cor_{r_{0}}} \EM{J}_2\cdot\EM{V}\mathrm d \mathbf x
		=\int_{\partial \Cor_{r_{0}}\cap \partial B_{r_0}}
		\EM{W}\cdot\pare{\nor\cros\EM{E}}\mathrm d \sigma
		+\int_{\partial \Cor_{r_{0}}\cap \partial B_{r_0}} \EM{V}\cdot\pare{\nor\cros\EM{H}}\mathrm d \sigma .
		\end{equation}	
		
		We first prove $\EM{J}_2(\mathbf 0)=\mathbf 0$. The conclusion for $\EM{J}_1(\mathbf 0)$ can be obtained similarly.  Since $\mathbf{J}_j\in C^\alpha(\overline{\Cor_{r_{0}}})^3$, we can write
		\begin{equation}\label{eq:f2}
		\EM{J}_2=\EM{J}_0+\tilde{\EM{J}},\quad \EM{J}_0=\EM{J}_2(\mathbf 0), 
		\end{equation}
		where $\tilde{\EM{J}}$ is a vector field satisfying 
		\[
		|\tilde{\EM{J}}(\mathbf x)|\le \|\EM{J}_2\|_{C^{\alpha }}|\mathbf x|^{\alpha},\quad \mathbf x\in \Cor_{r_{0}}.
		\]
		Substituting \eqref{eq:f2} into \eqref{eq:IntIdProof}, one has
		\begin{equation}\label{eq:TempProof}
		\begin{split}
		\int_{\Cor_{r_{0}}} \EM{J}_0\cdot \EM{V}\mathrm d \mathbf x
		=& 	\int_{\partial \Cor_{r_0}\cap \partial B_{r_0}} \EM{W}\cdot\pare{\nor\cros\EM{E}}\mathrm d \sigma
		+\int_{\partial \Cor_{r_0}\cap \partial B_{r_0}} \EM{V}\cdot\pare{\nor\cros\EM{H}}\mathrm d \sigma
		\\&-\int_{\Cor_{r_{0}}} \tilde{\EM{J}}\cdot\EM{V}\mathrm d \mathbf x
		-\int_{\Cor_{r_{0}}} \EM{J}_1\cdot\EM{W}\mathrm d \mathbf x.
		\end{split}
		\end{equation}
		
		Let $ \EM{V} $ and $ \EM{W} $ be defined by \eqref{eq:vw} which is a pair of solution to 
		the Maxwell system \eqref{eq:MaxwellTest}, where $  \boldsymbol{\rho} $ and $ \mathbf p $ are defined in \eqref{eq:r} satisfying \eqref{eq:innerdot}.
		
		Concerning the LHS of \eqref{eq:TempProof},  when $\tau$ is sufficiently large, from Lemma~\ref{lem:Nonvanish} we obtain that 
		\begin{equation}\label{eq:v}
		\abs{\int_{\Cor_{r_{0}}} \EM{J}_0 \cdot \EM{V}\mathrm d \mathbf x} = \abs{\EM{J}_0\cdot
			\mathbf p} \abs{\int_{\Cor_{r_{0}}} e^{\boldsymbol{\rho}\cdot \mathbf x} {\mathrm d} \mathbf x} \ge
		\abs{\EM{J}_0\cdot \mathbf p} C_{\Cor}(1+\frac{k^2}{\tau^2})^{-3/2} \tau^{-3}+\mathcal{O}(\tau^{-1}e^{-\frac{1}{2}r_{0}\tau \delta }). 
		\end{equation}
		 We shall show  that the RHS of \eqref{eq:TempProof} is bounded by $C\tau^{-(3+\alpha)}$ as $\tau \rightarrow +\infty $, where $C_{\Cor }$ is a positive constant and independent of  $ \tau $. 
		
		We first deal with the terms in \eqref{eq:TempProof} concerning $\EM{V}$. Since the apex  of $\Cor_{r_0}$ is the origin and the axis of $\Cor_{r_0}$ coincides with $ \mathbf x_{3}^{+} $, we can choose $\mathbf d=(0,0,-1)^\top$ and $\mathbf d^{\perp}=(\cos\varphi ,\sin \varphi , 0)^\top$,  where $\varphi\in (0,2\pi]$.  Hence  $\mathbf d$ fulfills the condition \eqref{eq:innerdot}  with $\delta =\cos \theta_0>0$ and $\theta_0 \in (0,\pi/2)$. 
		 Therefore using \eqref{eq:innerdot}  and \eqref{eq:f2} we have
		\begin{align}\label{eq:TempProof1}
		\abs{\int_{\Cor_{r_0}} \tilde{\EM{J}}\cdot\EM{V}\mathrm d \mathbf x}
		&\le \|\EM{J}_2\|_{C^\alpha}|\mathbf p|\int_{\Cor_{r_0}} |\mathbf x|^{\alpha}e^{\tau \mathbf d\cdot \mathbf x}{\mathrm d} \mathbf x
		\le 3\|\EM{J}_2\|_{C^\alpha}\tau^{-(3+\alpha)}\int_{\Cor}|\mathbf y|^{\alpha}e^{\mathbf d\cdot \mathbf y}{\mathrm d} \mathbf y \notag \\
		&\le 3\|\EM{J}_2\|_{C^\alpha}\tau^{-(3+\alpha)}\int_{\Cor}|\mathbf y|^{\alpha}e^{-\delta |\mathbf y|}{\mathrm d} \mathbf y\le C_{\Cor,\alpha}\|\EM{J}_2\|_{C^\alpha}\tau^{-(3+\alpha)},
		\end{align}
		where $C_{\Cor,\alpha}$ is a positive constant independently of $\tau$.

		
		For the boundary integral in \eqref{eq:TempProof}, by virtue  of  the trace theorem and Lemma \ref{lem:infty}  we have the estimate
		\begin{equation}\label{eq:TempProof2}
		\begin{split}
		\abs{\int_{\partial \Cor_{r_{0}}\cap \partial B_{r_0}} \EM{V}\cdot\pare{\nor\cros\EM{H}}\mathrm d  \sigma }
		&\le 3\int_{\partial \Cor_{r_{0}}\cap \partial B_{r_0}} \abs{\nor\cros\EM{H}} e^{-\delta \tau |\mathbf x|}{\mathrm d} \sigma
		\\&\le 3C_{\Cor,r_0}e^{-\delta r_0\tau}
		\|\EM{H}\|_{H(\curl,\, \Cor_{r_{0}})}.
		\end{split}	
		\end{equation}
		
		
Using  \eqref{eq:innerdot} and \eqref{eq:p norm}, 	 we have
		\begin{equation}\label{eq:w}
		\abs{\int_{\partial \Cor_{r_{0}}\cap \partial B_{r_0}} \EM{W}\cdot\pare{\nor\cros\EM{E}}\mathrm d \sigma}
		\le 3C_{\Cor,r_0}k^2\tau^{-1}e^{-\delta r_0\tau}
		\|\EM{E}\|_{H(\curl,\, \Cor_{r_{0}})}.
		\end{equation}
		
By using a similar argument as for deriving \eqref{eq:TempProof1}, one can obtain
		\begin{equation}\label{eq:ww}
		\begin{split}
		\abs{\int_{\Cor_{r_{0}}} \EM{J}_1\cdot\EM{W}\mathrm d \mathbf x}
		\le \omega\Ep_0\tau^{-1}\|\EM{J}_1\|_{C^{0}}
		\int_{\Cor_{r_{0}}} e^{\tau \mathbf d\cdot \mathbf x}{\mathrm d} \mathbf x
		\le C_{\Cor}\omega\Ep_0\|\EM{J}_1\|_{C^{0}}\tau^{-4}.
		\end{split}
		\end{equation}
		
		In view of \eqref{eq:v}, \eqref{eq:TempProof1}, \eqref{eq:TempProof2}, \eqref{eq:w} and \eqref{eq:ww}, and by virtue of \eqref{eq:TempProof}, one can show that
		\begin{equation}\label{eq:tau}
		\abs{\EM{J}_0\cdot \mathbf p} C_{\Cor}  \frac{\tau^{-3}}{1+\frac{k^2}{\tau^2}}+\mathcal{O}(\tau^{-1}e^{-\frac{1}{2}r_{0}\tau \delta })\le \tilde{C}_{\EM{J}_1,\EM{J}_2,\Cor_{r_0},\alpha}\tau^{-(3+\alpha)},
		\end{equation}
		where $ C_{\mathcal K} $ is a positive constant independently of $\tau$.
		
		Multiplying $ \tau^{3} $ on both side of \eqref{eq:tau}, let $ \tau \rightarrow +\infty $, we can deduce that
		\begin{equation}\label{eq:J_0}
		\lim_{\tau\to\infty} \EM{J}_0\cdot \mathbf p=0.
		\end{equation}
		
		Combining \eqref{eq:J_0} with Lemma \ref{lem:vani yang}, we can prove that  $\mathbf J_2(\mathbf 0)=\mathbf 0$.
%
 Finally, one can verify $\EM{J}_1(\mathbf 0)= 0$ in the same way  by taking
	\[
		\EM{V}(\mathbf x)=-\frac{1}{\im \omega\Ep_0} \boldsymbol{\rho}\cros \mathbf p e^{\boldsymbol{\rho}\cdot \mathbf x}\quad\text{and}\quad 
		\EM{W}(\mathbf x)=\mathbf pe^{\boldsymbol{\rho}\cdot \mathbf x}.
		\]
			
		The proof is complete. 
	\end{proof}
	
	In the following theorem we give a vanishing characterization  of non-radiating electric and magnetic sources associated with \eqref{eq:Maxwell1} near a conical corner. 
	
	\begin{thm}\label{thm:thm2.2}
	    Let $ \Omega \subset \mathbb R^3 $ be a bounded Lipschitz domain with a connected complement and $\mathbf x_0\in \partial \Omega$. Suppose that $ \Omega \cap B_{r_0}(\mathbf x_0)=\Cor_{\mathbf x_{0}}^{r_0} \cap B_{r_0}(\mathbf x_0) $ for some $ r_{0} \in \mathbb R_{+} $, where  $\Cor_{\mathbf x_{0}}^{r_0} $ is a conical corner defined by \eqref{eq:coner0}.  Let $ \mathbf{J}_1 $ and $ \mathbf{J}_2 $ be respectively  electric source and magnetic source both supported in $ \Omega $ fulfilling $  \mathbf{J}_1, \mathbf{J}_2 \in C^\alpha(\overline{ \Cor^{r_{0}}_{\mathbf x_0}} )^3 $. Consider the electromagnetic source scattering problem \eqref{eq:Maxwell1}.     
	     If $ \mathbf{J}_1 $ and $ \mathbf{J}_2 $ are non-radiating, namely $\mathbf{E}_\infty=\mathbf{H}_\infty\equiv\mathbf{0}$, then it holds that
	\begin{equation}\label{eq:cond2}
	\mathbf{J}_1(\mathbf x_0)=\mathbf{J}_2(\mathbf x_0)=\mathbf  0. 
	\end{equation}
    \end{thm}

\begin{proof}
	Let  $ (\mathbf E, \mathbf H)\in H({\rm curl},\Omega)  \times H({\rm curl},\Omega ) $ be   pair of solutions for  the Maxwell system \eqref{eq:MaxwellO} associated with electric and magnetic sources $\mathbf J_1$ and   $\mathbf J_2$. Since $ \mathbf{J}_1 $ and $ \mathbf{J}_2 $ are radiationless, one has the far-field pattern $ (\mathbf E_{\infty},\mathbf H_{\infty}) \equiv \mathbf 0 $, and then from Rellich's theorem (cf. \cite{CoK13}), we can obtain that $ \mathbf E=\mathbf H=\mathbf 0 $ in $ \mathbb R^3\backslash \bar{\Omega} $, which imply \eqref{eq:MaxwellCor0}.   Finally, by using Theorem \ref{thm:main0}, we can prove this theorem. 
\end{proof}


As a direct consequence of Theorem~\ref{thm:thm2.2}, we can derive the following result which shows that under generic conditions, a pair of electromagnetic sources must radiate a nontrivial scattering pattern if they possess a conical singularity on their supports, namely they must be visible. 


\begin{thm}\label{thm:41}
	Consider the electromagnetic source scattering problems \eqref{eq:Maxwell1}. If the support of an electric source ${\mathbf J}_1$ or a magnetic source ${\mathbf J}_2$  has a conical corner $\mathcal K_{\mathbf x_0}^{r_0}$ as described in \eqref{eq:coner0},   where either ${\mathbf J}_1\in C^{\alpha} ( \overline {\mathcal K_{\mathbf x_0}^{r_0}}  )$ with $ {\mathbf J}_1 (\mathbf x_0) \neq \mathbf 0 $  or  ${\mathbf J}_2\in C^{\alpha} (\overline{ \mathcal K_{\mathbf x_0}^{r_0} } )$ with ${\mathbf J}_2 (\mathbf x_0) \neq \mathbf 0 $ is fulfilled, then ${\mathbf J}_1$ or ${\mathbf J}_2$ always radiate a nontrivial scattering pattern for any angular frequency $\omega$.

	\end{thm}

\subsection{Inverse source scattering}

In the following  we shall deal with the unique recovery results on the nonlinear inverse scattering  \eqref{eq:ip2}. Before that we introduce the admissible electromagnetic source configurations.  

\begin{defn}\label{def:adm source}   
	Suppose that $\Omega$ is a bounded Lipschitz domain with a connected complement and $ \mathbf J_{1},\ \mathbf J_{2} \in L_{\rm loc}^2(\mathbb{R}^3; \mathbb C^3) $ with $ {\rm supp}(\mathbf J_{1})=\Omega ,\ {\rm supp}(\mathbf J_{2})=\Omega $. If $\Omega$ has a conical corner $\Cor_{\mathbf x_0}^{r_0}$ described by \eqref{eq:coner0} such that $ \Omega \cap B_{r_0}(\mathbf x_0)=\Cor_{\mathbf x_{0}}^{r_0} \cap B_{r_0}(\mathbf x_0) $ for some $ r_{0} \in \mathbb R_{+} $, where $ \mathbf J_{1}$ and $\mathbf J_{2} $ have the H\"older continuous  regularity near the underlying conical corner fulfilling the condition
	\begin{align}\label{eq:JJ}
		\mbox{$\mathbf{J}_1(\mathbf x_0)\neq \mathbf 0$ or $\mathbf{J}_2(\mathbf x_0)\neq  \mathbf 0$. } 
	\end{align}
	 Then $ (\Omega; \mathbf{J}_1, \mathbf{J}_2) $ is said to belong to the set of admissible electromagnetic source configurations.
\end{defn}

A local unique recovery result concerning the inverse source shape  \eqref{eq:ip2} by a single far field measurement can be established by using Theorem~\ref{thm:main0}.

\begin{thm}\label{thm:main03}
	Suppose that $(\Omega; \mathbf{J}_1, \mathbf{J}_2)$ and $(\Omega'; \mathbf{J}'_1, \mathbf{J}'_2)$ are two admissible  electromagnetic source configurations described in Definition~\ref{def:adm source}, where $ \mathbf J_{1},\ \mathbf J_{2} $ are supported in $\Omega$ and $\mathbf{J}'_1,\mathbf{J}'_2$ in $\Omega^{'}$.   
	Let $\mathbf{E}_\infty$ and $\mathbf{E}_\infty'$ be the electric far-field patterns associated with $(\Omega; \mathbf{J}_1, \mathbf{J}_2)$ and $(\Omega'; \mathbf{J}'_1, \mathbf{J}'_2)$  respectively. If 
	\begin{equation}\label{eq:EE}
		\mathbf{E}_\infty(\hat {\mathbf x})=\mathbf{E}_\infty'(\hat {\mathbf x}) \quad \text{for all}\quad \hat {\mathbf x}\in\mathbb{S}^2,
	\end{equation}
	 then the set difference
	\begin{equation}\label{eq:diff01}
	\Omega\Delta\Omega':=(\Omega\setminus\overline{\Omega'})\cup (\Omega'\setminus\overline{\Omega}), 
	\end{equation}
	cannot contain a conical corner.	
\end{thm}

\begin{proof}
	We prove the theorem by contradiction. Suppose that  $(\mathbf{E}, \mathbf{H})$ and $(\mathbf{E}',\mathbf{H}')$ are the electromagnetic fields of the Maxwell system \eqref{eq:Maxwell1} associated with $(\Omega; \mathbf{J}_1, \mathbf{J}_2)$ and $(\Omega'; \mathbf{J}'_1, \mathbf{J}'_2)$  respectively. Let 
	\begin{equation}\label{eq:G}
		\mathbf{G}:=\mathbb{R}^3\setminus\overline{\Omega\cup\Omega'}.
	\end{equation}
	By contradiction, without loss of generality, we suppose that  $r_0\in\mathbb{R}_+$ are sufficient small such that the conical corner $\mathcal{K}_{\mathbf x_0}^{r_0}\subset\overline{\Omega}\setminus\Omega'$ 
	, where $ \mathcal{K}_{\mathbf x_0}^{r_0}=\mathcal{K}_{r_0} $ is defined by \eqref{eq:coner0} and $\mathbf x_0\in\overline{\Omega}\setminus\Omega'$ satisfying  \eqref{eq:JJ}. According to \eqref{eq:EE}, by the Rellich theorem, we readily have that
	\begin{equation}\label{eq:ar1}
	(\mathbf{E}, \mathbf{H})=(\mathbf{E}', \mathbf{H}')\quad\mbox{in}\ \ \mathbf{G}. 
	\end{equation}
	 Therefore, we can obtain that
	\begin{equation}\label{eq:ar2}
	\begin{cases}
	&\nabla\wedge \EM{E}(\mathbf x)-\im\omega \Mp_0 \EM{H}(\mathbf x)=\EM{J}_1(\mathbf x),\quad\, \mathbf x\in\mathcal{K}_{\mathbf x_0}^{r_0},\\
	& \nabla\wedge \EM{H}(\mathbf x)+\im\omega\Ep_0 \EM{E}(\mathbf x)=\EM{J}_2(\mathbf x),\quad\ \mathbf x\in\mathcal{K}_{\mathbf x_0}^{r_0},
	\end{cases}
	\end{equation}
	and 
	\begin{equation}\label{eq:ar3}
	\begin{cases}
	&\nabla\wedge \EM{E}'(\mathbf x)-\im\omega \Mp_0 \EM{H}'(\mathbf x)=\mathbf 0,\quad\, \mathbf x\in\mathcal{K}_{\mathbf x_0}^{r_0},\\
	& \nabla\wedge \EM{H}'(\mathbf x)+\im\omega\Ep_0 \EM{E}'(\mathbf x)=\mathbf 0,\quad\ \mathbf x\in\mathcal{K}_{\mathbf x_0}^{r_0},
	\end{cases}
	\end{equation}\label{eq:EEHH}
Let 
	\begin{align}\label{eq:E tilde E}
			\widetilde{\EM{E}}=\mathbf{E}-\mathbf{E}',\quad\widetilde{\EM{H}}=\mathbf{H}-\mathbf{H}'.
	\end{align}
	By \eqref{eq:ar1}--\eqref{eq:ar3}, it readily holds that 
	\begin{equation}\label{eq:ar4}
	\begin{cases}
	&\nabla\wedge \widetilde{\EM{E}}(\mathbf x)-\im\omega \Mp_0 \widetilde{\EM{H}}(\mathbf x)=\EM{J}_1(\mathbf x),\quad\, \mathbf x\in\mathcal{K}_{\mathbf x_0}^{r_0},\\
	& \nabla\wedge \widetilde{\EM{H}}(\mathbf x)+\im\omega\Ep_0 \widetilde{\EM{E}}(\mathbf x)=\EM{J}_2(\mathbf x),\quad\ \mathbf x\in\mathcal{K}_{\mathbf x_0}^{r_0},\\
	&\nu(\mathbf x)\wedge\widetilde{\mathbf{E}}(\mathbf x)=\nu(\mathbf x)\wedge\widetilde{\mathbf{H}}(\mathbf x)=\mathbf 0,\quad \mathbf x\in\partial\mathcal{K}_{\mathbf x_0}^{r_0}\backslash \partial B_{r_0}(\mathbf x_0). 
	\end{cases}
	\end{equation}
	Therefore, by Theorem~\ref{thm:main0}, we can derive that
	\[
	\mathbf{J}_1(\mathbf x_0)=\mathbf{J}_2(\mathbf x_0)=\mathbf 0,
	\]
	which is a contradiction to \eqref{eq:JJ}. 
\end{proof}

In Theorem \ref{thm:global source}, we shall establish a global unique identifiability results for inverse electromagnetic source problem with certain a-prior knowledge  on the sources. Before that, we introduce the definition of  an admissible electromagnetic source configuration of coronal shape.

\begin{defn}\label{def:finite c}
    Suppose that $ D $ is a convex bounded Lipschitz domain with a connect complement $ \mathbb R^3\backslash \bar{D}$. If there exist finite many strictly convex conical cone $ \mathcal K_{\mathbf x_j,\theta_{j}}(j=1,...,\ell,\ell \in \mathbb N) $ defined by \eqref{eq:ccone} ($2\theta_j$ is the opening angle of the cone $\mathcal K_{\mathbf x_j,\theta_{j}}$) satisfying
    \begin{enumerate}
    	\item the apex $ \mathbf x_j $ of $ \mathcal K_{\mathbf x_j,\theta_{j}} $ satisfies $ \mathbf x_j\in \mathbb R^3\backslash \bar{D} $ and let $ \mathcal K_{\mathbf x_j,\theta_{j}}\backslash \bar{D} $ be denoted by $ \mathcal K_{\mathbf x_j,\theta_{j}}^{(1)} $, where the apex $ \mathbf x_j $ belongs to the strictly convex bounded conical cone $ \mathcal K_{\mathbf x_j,\theta_{j}}^{(1)} $,
    	\item $ \partial \overline{\mathcal K_{\mathbf x_j,\theta_{j}}^{(1)}} \backslash \partial \overline{\mathcal K_{\mathbf x_j,\theta_{j}}}\subseteq
    	 \partial \bar{D}  $ and $ \cap_{j=1}^{\ell}\partial \overline{\mathcal K_{\mathbf x_j,\theta_{j}}^{(1)}} \backslash \partial \overline{\mathcal K_{\mathbf x_j,\theta_{j}}}=\emptyset $, 
    \end{enumerate}
then $ \Omega :=\cup_{j=1}^{\ell} \mathcal K_{\mathbf x_j,\theta_{j}} \cup D $ is said to be bounded Lipschitz domain of coronal shape with the corresponding  conical corners $ \mathcal K_{\mathbf x_j,\theta_{j}}^{(1)}\ (j=1,...,\ell) $; see Fig.~\ref{fig:coordinate1} for a schematic illustration. Let $ \mathbf{J_{1}}, \mathbf{J_{2}} \in L^2_{\rm loc}(\mathbb R^3) $ with $ {\rm supp}(\mathbf J_{1})=\Omega $ and $ {\rm supp}(\mathbf J_{2})=\Omega $, where $ \mathbf J_{1}$ and $\mathbf J_{2} $ have H\"older continuous regularity in $\mathcal K_{\mathbf x_j,\theta_{j}}^{(1)}\ (j=1,...,\ell)$ fulfilling
 \begin{equation}\label{eq:JJ1}
	\begin{split}
	&\mathbf J_{1}({\mathbf x_{j}})\not=\mathbf 0 \quad \mbox{or}\quad \mathbf J_{2}({\mathbf x_{j}})\not=\mathbf 0,\quad \forall j=1,...,\ell. 
	\end{split}	
	\end{equation} 
Then $ (\Omega ;\mathbf J_1,\mathbf J_2) $ is said to belong to the admissible electromagnetic source configuration of coronal shape. 
\end{defn}

\begin{thm}\label{thm:global source}
	Suppose that $ (\Omega ;\mathbf J_1,\mathbf J_2) $ and $ (\Omega' ;\mathbf J_1^{'},\mathbf J_2^{'}) $ are two admissible electromagnetic source configurations of coronal shape described in Definition \ref{def:finite c}, where $\Omega = \cup_{j=1}^{\ell} \mathcal K_{\mathbf x_j,\theta_{j}} \cup D $, $ \Omega^{'} = \cup_{j=1}^{\ell^{'}} \mathcal K_{\mathbf x_j^{'},\theta_{j}^{'}} \cup D^{'} $ and $ \mathbf J_i \in C^\alpha(\mathcal K_{\mathbf x_j, \mathbf \theta_{j}}) $. 		Let
	$ (\mathbf{E}_\infty, \mathbf{H}_\infty) $ and $ (\mathbf{E}_\infty^{'}, \mathbf{H}_\infty^{'}) $ be the far field patterns associated with $ (\Omega ;\mathbf J_1,\mathbf J_2) $ and $ (\Omega ;\mathbf J_1^{'},\mathbf J_2^{'}) $ respectively. If 
	\begin{equation}\label{eq:EH}
		(\mathbf{E}_\infty, \mathbf{H}_\infty)=(\mathbf{E}_\infty^{'}, \mathbf{H}_\infty^{'})\quad \text{for all} \quad \hat{\mathbf x} \in \mathbb S^2,\quad D=D',
	\end{equation}
	and \begin{equation}\label{eq:J1J2}
		\mathbf J_{1}({\mathbf x_{j}})\not=\mathbf J^{'}_{1}({\mathbf x^{'}_{i}})\quad \mbox{or} \quad \mathbf J_{2}({\mathbf x_{j}})\not=\mathbf J^{'}_{2}({\mathbf x^{'}_{i}}),
	\end{equation}
	where $ \mathbf x_{j}=\mathbf x^{'}_{i} $ for some $ j\in \{1,...,\ell\} $ and $ i\in \{1,...,\ell^{'}\} $, 
	then $\Omega=\Omega'.$ 
	
	If \eqref{eq:EH} is satisfied and 
	\begin{align}\label{eq:source 239}
		\mbox{ $\theta_j=\theta_i'$ for  any $j\in \{1,\ldots,\ell\}$ and $i\in \{1,\ldots, \ell'\}$ fulfilling $\mathbf x_j=\mathbf x_i'$,}
	\end{align}
	then $\Omega=\Omega'$, 
	\begin{align}\label{eq:conclusion 240}
		\mbox{$\mathbf J_{1}({\mathbf x_{j}}) =\mathbf J^{'}_{1}({\mathbf x_{j}})$\quad and \quad $\mathbf J_{2}({\mathbf x_{j}})=\mathbf J^{'}_{2}({\mathbf x_{j}})$, }
	\end{align}
	where $j\in \{1,\ldots, \ell\}$. 
	
\end{thm}
\begin{proof}
	We prove this theorem by the contradiction. Recall that $ D=D' $. Suppose that $ \ell\not=\ell^{'} $ or $ \mathbf x_{j}\not=\mathbf x^{'}_{j} $, we can see that $ \Omega \Delta \Omega' $ defined in \eqref{eq:diff01} has a conical corner, where the underlying source is non-vanishing at the corresponding corner by virtue of \eqref{eq:JJ1}, which is contradict to   \eqref{eq:JJ}. According to Theorem~\ref{thm:main03}, we have $ \ell=\ell^{'} $ and $ \mathbf x_{j}=\mathbf x^{'}_{j}\ (j=1,...,\ell) $. 
	
	In the following we prove that    $ \theta_{j}=\theta_{j}^{'} ,\ \forall  j\in \{1,...,\ell\} $. By contradiction,  there exists an index $ j_0 \in \{1,...,\ell\} $ such that $ \theta_{j_{0}} \neq \theta_{j_{0}}^{'} $. Without loss of generality, we may suppose that $ \theta_{j_{0}}^{'} <\theta_{j_{0}} $. By virtue of \eqref{eq:EH}, we can obtain that $\mathbf  E=\mathbf E' $ and $ \mathbf H=\mathbf H' $ in $ \mathbf{G} $, where $ \mathbf{G} $ is defined by \eqref{eq:G}. Therefore we have 
	\begin{equation}\label{eq:J}
		\begin{cases}
		&\nabla\wedge \EM{E}(\mathbf x)-\im\omega \Mp_0 \EM{H}(\mathbf x)=\EM{J}_1(\mathbf x),\quad\, \mathbf x\in\mathcal{K}_{\mathbf x_{j_{0}},\theta_{j}}\cap B_{r_{0}}(\mathbf x_{j_{0}}),\\
		& \nabla\wedge \EM{H}(\mathbf x)+\im\omega\Ep_0 \EM{E}(\mathbf x)=\EM{J}_2(\mathbf x),\quad\ \mathbf x\in\mathcal{K}_{\mathbf x_{j_{0}},\theta_{j}}\cap B_{r_{0}}(\mathbf x_{j_{0}}),
		\end{cases}
	\end{equation}
and 
\begin{equation}\label{eq:J'}
	\begin{cases}
	&\nabla\wedge \EM{E}'(\mathbf x)-\im\omega \Mp_0 \EM{H}'(\mathbf x)=\chi_{\Omega'} \mathbf J_{1}^{'}(\mathbf x) ,\quad\, \mathbf x\in\mathcal{K}_{\mathbf x_{j_{0}},\theta_{j}}\cap B_{r_{0}}(\mathbf x_{j_{0}}),\\
	& \nabla\wedge \EM{H}'(\mathbf x)+\im\omega\Ep_0 \EM{E}'(\mathbf x)=\chi_{\Omega'} \mathbf J_{2}^{'}(\mathbf x),\quad\ \mathbf x\in\mathcal{K}_{\mathbf x_{j_{0}},\theta_{j}}\cap B_{r_{0}}(\mathbf x_{j_{0}}),
	\end{cases}
\end{equation}
where $ r_0 \in \mathbb R_+ $ is suffiently small such that $ \Cor_{\mathbf x_{j_{0}},\theta_{j}} \cap B_{h}(\mathbf x_{j_{0}}) \subset \mathbf R^3\backslash \bar{D}$. By virtue of \eqref{eq:J} and \eqref{eq:J'}, we can obtain that 
\begin{equation}
	\begin{cases}
	&\nabla\wedge \widetilde{\EM{E}}(\mathbf x)-\im\omega \Mp_0 \widetilde{\EM{H}}(\mathbf x)=\EM{J}_1(\mathbf x)-\chi_{\Omega'} \mathbf J_{1}^{'}(\mathbf x),\quad\, \mathbf x\in\mathcal{K}_{\mathbf x_{j_{0}},\theta_{j}}\cap B_{r_{0}}(\mathbf x_{j_{0}}),\\
	& \nabla\wedge \widetilde{\EM{H}}(\mathbf x)+\im\omega\Ep_0 \widetilde{\EM{E}}(\mathbf x)=\EM{J}_2(\mathbf x)-\chi_{\Omega'} \mathbf J_{2}^{'}(\mathbf x),\quad\mathbf  x\in\mathcal{K}_{\mathbf x_{j_{0}},\theta_{j}}\cap B_{r_{0}}(\mathbf x_{j_{0}}),\\
	&\nu(\mathbf x)\wedge\widetilde{\mathbf{E}}(\mathbf x)=\nu(\mathbf x)\wedge\widetilde{\mathbf{H}}(\mathbf x)=\mathbf 0,\quad \mathbf x\in \partial(\mathcal{K}_{\mathbf x_{j_{0}},\theta_{j}}\cap B_{r_{0}}(\mathbf x_{j_{0}}))\backslash \partial B_{r_{0}}(\mathbf x_{j_{0}}).
	\end{cases}
\end{equation}
where $ \widetilde{\EM{E}} $ and $ \widetilde{\EM{H}} $ are defined in \eqref{eq:E tilde E}.  Using  Theorem~\ref{thm:main0}, one has
\begin{equation}\notag
	\mathbf J_{1}({\mathbf x_{j_0}})=\mathbf J^{'}_{1}({\mathbf x_{j_0}}),\ \mathbf J_{2}({{\mathbf x}_{j_0}})=\mathbf J^{'}_{2}({\mathbf x_{j_0}}),
\end{equation}
which is a contradiction to \eqref{eq:J1J2}.

When \eqref{eq:source 239}  is fulfilled, one can readily  obtain \eqref{eq:conclusion 240}  by using Theorem \ref{thm:main0}. 
\end{proof} 

As discussed in the introduction, by Theorems~\ref{thm:41} and \ref{thm:global source}, one can conclude that if an admissible electromagnetic source of coronal shape is visible, then it can be uniquely recovered by using the visible far-field pattern.

\section{Transparent/invisible scatterers and inverse medium scattering}\label{sec:3}

Consider the electromagnetic medium scattering problem  formulated by \eqref{eq:Maxwell2}, where the physical parameters  of the medium scatterer $(\Omega; \varepsilon,\mu,\sigma)$ is described by \eqref{eq:Maxwell1}. In this section, we first establish a geometrical charaterisation of transparent/invisible mediums near a conical corner when the  total wave field holds a H\"older regularity condition near the corresponding  conical corner. When the medium scatterer $\Omega$ is simply connected and transparent/invisible, namely $\mathbf{E}_\infty\equiv \mathbf 0$ or $\mathbf{H}_\infty\equiv \mathbf 0$, by Rellich's theorem, one can directly know that  $(\mathbf{E}, \mathbf{H})=\mathbf 0$ in $\mathbb{R}^3\setminus\overline{\Omega}$, where  $(\mathbf{E}, \mathbf{H})$ is the scattered wave field of \eqref{eq:Maxwell2}. Hence, one can show that the total and incident wave fields, namely  $(\mathbf{E}^t, \mathbf{H}^t)$ and $(\mathbf{E}^i, \mathbf{H}^i)$, satisfy the following system:  
\begin{equation}\label{eq:ITP old}
		\begin{cases}
			\nabla\wedge \EM{E}^t-\mathrm{i}\omega\Mp \EM{H}^t=\mathbf 0,\quad
			\nabla \wedge \EM{H}^t+\mathrm{i}\omega\Eci \EM{E}^t=\mathbf 0
			&\mbox{in\ \ $\Omega$},\\
			\nabla\wedge \EM{E}^i-\mathrm{i}\omega\Mp_0 \EM{H}^i=\mathbf 0,\quad
			\nabla\wedge \EM{H}^i+\mathrm{i}\omega\Ep_0 \EM{E}^i=\mathbf 0
			&\mbox{in\ \ $\Omega$},\\
			\quad\nu\cros \EM{E}^t=\nu\cros \EM{E}^i,\quad
			\qquad\nu\cros \EM{H}^t=\nu\cros \EM{H}^i
			&\mbox{on\ \ $\partial \Omega$},
		\end{cases}
	\end{equation}
where $\gamma$ is defined in \eqref{eq:Maxwell2}. \eqref{eq:ITP old} is known as the electromagnetic transmission eigenvalue problem in the literature and we shall present more discussion in Section~\ref{sec:4} in what follows.  


\subsection{Geometrical characterisation of transparent/invisible mediums}

We first establish a geometrical characterisation of transparent/invisible mediums near the a conical corner under certain regularity conditions. 




\begin{thm}\label{thm:main5}
Consider the electromagnetic medium  scattering problem \eqref{eq:Maxwell2} with the associated 	electromagnetic medium  scatterer $ (\Omega ; \varepsilon,\mu,\sigma ) $,  where $ \Omega \subset \mathbb R^3 $ be a bounded Lipschitz domain with a connected complement and $\mathbf x_0\in \partial \Omega$. Suppose that $ \Omega \cap B_{r_0}(\mathbf x_0)=\Cor_{\mathbf x_{0}}^{r_0} \cap B_{r_0}(\mathbf x_0) $ for some $ r_{0} \in \mathbb R_{+} $, where  $\Cor_{\mathbf x_{0}}^{r_0} $ is a conical corner defined by \eqref{eq:ccone}. Assume that
\begin{equation}\label{eq:rega1}
	(\mu(\mathbf x)-\mu_0)\mathbf{H}^t,\ (\gamma(\mathbf x)-\varepsilon_0)\mathbf{E}^t\in C^\alpha(\overline{\mathcal{K}^{r_0}_{\mathbf x_0}})^3,
	\end{equation} 
	for some $\alpha\in (0, 1)$, where $(\mathbf{E}^t, \mathbf{H}^t)$ is  the total wave field of \eqref{eq:Maxwell2}.  If  $ (\Omega ; \varepsilon,\mu,\sigma ) $ is transparent/invisible, namely, $ (\mathbf E_{\infty},\mathbf H_{\infty}) \equiv \mathbf 0 $, 
	then there holds
	\begin{equation}\label{eq:rega2}
	(\mu(\mathbf x_0)-\mu_0)\mathbf{H}^t(\mathbf x_0)=(\gamma(\mathbf x_0)-\varepsilon_0)\mathbf{E}^t(\mathbf x_0)=\mathbf 0. 
	\end{equation}
	Furthermore, if $\mu(\mathbf x_0)\neq \mu_0\neq 0$ and $\gamma(\mathbf x_0)\neq \varepsilon_0$, one has $\mathbf{H}^t(\mathbf x_0)=\mathbf{E}^t(\mathbf x_0)=\mathbf 0.$
\end{thm}

\begin{proof}
Since $(\mathbf E_{\infty},\mathbf H_{\infty}) \equiv \mathbf 0$, by Rellich theorem, one can directly know that  $(\mathbf{E}, \mathbf{H})=\mathbf 0$ in $\mathbb{R}^3\setminus\overline{\Omega}$. Hence, it arrives at  \eqref{eq:ITP old}. According  to \eqref{eq:ITP old}, by straightforward calculations, one can show   that $(\widehat{\mathbf{E}},\widehat{\mathbf{H}}):=(\mathbf{E}^t,\mathbf{H}^t)-(\mathbf{E}^i,\mathbf{H}^i)$ satisfies 
	\begin{equation}\label{eq:rega3}
	\left\{
	\begin{split}
	& \nabla\wedge  \widehat{\EM{E}}-\im\omega \Mp_0 \widehat{\EM{H}}=\EM{J}_{1} \quad \mbox{in $\Cor^{r_0}_{\mathbf x_0}$},\\
	& \, \nabla\wedge \widehat{\EM{H}}+\im\omega\Ep_0 \widehat{\EM{E}}=\EM{J}_{2} \quad \mbox{in $\Cor^{r_0}_{\mathbf x_0}$},\\
	&\, \nor\cros \widehat{\EM{E}}=\nor\cros\widehat{\EM{H}}=\mathbf 0 \hspace*{0.8cm}\mbox{on $\partial \Cor^{r_0}_{\mathbf x_0}\setminus \partial B_{r_0}(\mathbf x_0)$ ,} 
	\end{split}\right.
	\end{equation}
	where
	\begin{equation}\label{eq:rega4}
	\mathbf{J}_1= \mathrm i \omega (\mu(\mathbf x)-\mu_0)\mathbf{H}^t\quad\mbox{and}\quad\mathbf{J}_2=\mathrm i \omega (\gamma(\mathbf x)-\varepsilon_0)\mathbf{E}^t. 
	\end{equation}
	Hence, under the assumption  \eqref{eq:rega1}, by Theorem~\ref{thm:main0}, one readily has \eqref{eq:rega2}. 
\end{proof}

Similar to Theorem \ref{thm:41}, where the electromagnetic source with a conical corner always radiates a nontrivial far-field pattern, we shall reveal that an electromagnetic medium scatter containing a conical corner scatters any incident wave nontrivially, namely it must be visible. In proving such a result, we shall first need to establish a certain regularity property of the scattering problem, which will be given in the next subsection. Hence, we postpone this result to the end of the next subsection and present it in Theorem \ref{thm:medium conical scat} in what follows.


\subsection{Inverse medium scattering}\label{subsec:in me}

In this subsection, we establish several unique identifiability results for  an admissible electromagnetic medium scatterer by a single far field measurement under certain physical assumptions.

We first introduce the admissible class for the electromagnetic medium scatterers.

\begin{defn}\label{def:adm}
	Let $ (\Omega ; \varepsilon,\mu,\sigma ) $ be an electromagnetic medium   scatterer associate with \eqref{eq:Maxwell2}. Denote $\gamma(\mathbf x)=\varepsilon(\mathbf x)+\mathrm{i}\sigma(\mathbf x)/\omega$. Consider the electromagnetic medium scattering  \eqref{eq:Maxwell2} and $ (\mathbf E^t,\mathbf H^t) $ is the total wave field therin. The scatterer is said to be admissible if it fulfills the following conditions:
	\begin{enumerate}
		\item $ \Omega $ is a bounded simply connected Lipschitz domain in $ \mathbb R^3 $. The electric permittivity  $\varepsilon$, magnetic permeability $\mu$ and electric conductivity $\sigma$ associated with  the medium scatterer $\Omega$ satisfy the following condition 
		\begin{align}\label{eq:adm para}
		\begin{split}
			&\varepsilon\in L^\infty(\Omega; \mathbb{R}_+)\cap H^1(\Omega; \mathbb{R}_+),  \quad  \mu\in L^\infty(\Omega; \mathbb{R}_+)\cap H^1(\Omega; \mathbb{R}_+),\\
			&\sigma\in L^\infty(\Omega; \mathbb{R}_+^0)\cap H^1(\Omega; \mathbb{R}_+^0).
				\end{split}
		\end{align}
		\item If $\Omega$ has a conical corner $ 	\Cor^{r_0}_{\mathbf x_0} \Subset \Omega $ with the form \eqref{eq:coner0},  where $ \mathbf x_0\in  \partial\Omega $ is the apex of the underlying corner $ \Cor^{r_{0}}_{\mathbf x_0} $  with a sufficient small $ r_{0} \in \mathbb R_{+} $, then $ \mu$ and $ \gamma $ fulfill the following condition
		\begin{equation}
		\mu(\mathbf  x )=\mu_1, \gamma(\mathbf x)=\gamma_1,\quad \forall \mathbf x\in  	\Cor^{r_0}_{\mathbf x_0}  , 
		\end{equation}
		where $\mu_1$ and $\gamma_1$ are positive constants satisfying $ \mu_1\not=\mu_0$ and $\gamma_1\not=\varepsilon_0 $. 

		\item Either $ \mathbf E^t $ or $\mathbf H^t$ is non-vanishing everywhere in the sense that for any $ \mathbf x\in \mathbb R^3 $,
		\begin{equation}\label{eq:nn}
		\begin{split}
		&\lim_{\rho \to +0}\frac{1}{m(B(\mathbf x,\rho))}\int_{ B(\mathbf x,\rho) } | \mathbf 
		G (\mathbf x) | \mathrm d \mathbf x \not= 0,\quad \mathbf G=\mathbf E^t \quad \mbox{or} \quad \mathbf G=\mathbf H^t. 
		\end{split}
		\end{equation}
	\end{enumerate}
\end{defn}

\begin{rem}
The assumption (3) in Definition \ref{def:adm} can be fulfilled in certain physical scenario. For example, when $\omega \cdot {\rm diam}(\Omega ) \ll 1 $,  the diameter of the scatterer is far smaller than the incident wavelength, which implies that the scattered wave field is not dominant compared with the incident wave. Nevertheless, we shall not investigate under what more general physical applications the assumption (3) may be satisfied in this paper. 
\end{rem}

\begin{rem}
	We emphasize that \eqref{eq:adm para} plays an important role in proving the H\"older continuous  regularity of the total  wave field at a conical corner point in Lemma \ref{lem:36 max}.  However, the assumption \eqref{eq:adm para} for the physical parameters $\varepsilon$, $\mu$ and $\sigma$ can be replaced with that they are piecewise constants in $\Omega$, which implies that Lemma \ref{lem:36 max} is valid under this situation. Hence the unique identifiability for the medium shape determination  by a single measurement in our subsequent discussions hold for the case that  $\varepsilon$, $\mu$ and $\sigma$ are piecewise constants in $\Omega$. 
\end{rem}

In the following theorem, we establish a local unique recovery on the shape determination for an admissible scatterer by a single far field measurement. Before that we first show local regularity results on the solutions to \eqref{eq:Maxwell2} in Lemma \ref{lem:36 max}, where the medium scatter is admissible. 

\begin{lem}\cite[Theorems 1 and 2]{Alb}\label{lem:31 regu}
	Let  $\Omega$ be a bounded and connected open set in $\mathbb R^3$, with $C^{1,1}$ boundary. Let $\varepsilon,\mu\in L^\infty(\Omega; \mathbb C^{3 \times 3} ) $ be two bounded complex matrix-valued functions with uniformly positive definite real parts and symmetric imaginary parts. For a given frequency $\omega \in \mathbb C \backslash \{0\} $ and current sources $\mathbf J_e$ and $\mathbf J_m$ in $L^2(\Omega; \mathbb C^3)$, let $(\widehat{\mathbf E},\widehat{\mathbf H})$ in $H({\rm curl},\Omega)$ be the weak solution to the following  time-harmonic anisotropic Maxwell's equations
	\begin{equation}\label{eq:regu lem}
	\left\{
	\begin{split}
	& \nabla\wedge  \widehat{\EM{E}}-\im\omega \Mp\widehat{\EM{H}}=\EM{J}_{m} \quad \mbox{in $\Omega$},\\
	& \, \nabla\wedge \widehat{\EM{H}}+\im\omega\Ep \widehat{\EM{E}}=\EM{J}_{e} \hspace{0.5cm} \mbox{in $\Omega$},\\
	&\, \nor\cros \widehat{\EM{E}}=\nor\cros\widehat{\EM{G}} \hspace*{1.5cm}\mbox{on $\partial  \Omega $, } 
	\end{split}\right.
	\end{equation}
	where $\widehat{ \mathbf G}\in H({\rm  curl},\Omega )$.

	If $\varepsilon \in W^{1,3+\delta}(\Omega;\mathbb C^{3\times  3} )$ and the source terms $\mathbf J_m, \mathbf J_e$, and $\widehat{\mathbf G} $ satisfy 
	\begin{align}
		\mathbf J_m \in L^p(\Omega;\mathbb C^3 ),\quad \mathbf J_e\in W^{1,p}({\rm div},\Omega ),\quad{and}\quad \widehat{\mathbf G}\in W^{1,p}(\Omega; \mathbb C^3 )
	\end{align}
	for some $p\geq 2$, where $W^{N,p}({\rm div},\Omega )=\{\mathbf V\in W^{N-1,p}(\Omega;\mathbb C^3 ):\nabla \cdot \mathbf V \in  W^{N-1,p}(\Omega;\mathbb C )\}$, then $\widehat{\mathbf E}\in H^1(\Omega;\mathbb C^3 ) $. 
	
	If $\mu \in W^{1,3+\delta}(\Omega;\mathbb C^{3\times  3} )$ and the source terms $\mathbf J_m, \mathbf J_e$, and $\widehat{\mathbf G} $ satisfy 
	\begin{align}
	\begin{split}
		\mathbf J_e \in L^p(\Omega;\mathbb C^3 ),\ \mathbf J_m \in W^{1,p}({\rm div},\Omega ),\ & \mathbf J_m \cdot\nu \in W^{1-\frac{1}{p}}(\partial  \Omega;\mathbb C ),   \\ &\mbox{and}\quad \widehat{\mathbf G}\in W^{1,p}(\Omega; \mathbb C^3 )
			\end{split}
	\end{align}
	for some $p\geq 2$, then $\widehat{\mathbf H}\in H^1(\Omega;\mathbb C^3 ) $.

\end{lem}

\begin{lem}\label{lem:36 max}
Consider the electromagnetic scattering problem \eqref{eq:Maxwell2}.   Let 
\begin{align}\label{eq:EH curl cond}
	 (\mathbf E^t,\mathbf H^t) \in  H_{\mathrm{loc}}(\mathrm{curl}, \mathbb{R}^3)\times H_{\mathrm{loc}}(\mathrm{curl}, \mathbb{R}^3)
\end{align}
 be the total wave field associated with the admissible scatterer $ (\Omega ; \varepsilon,\mu,\sigma )$. Assume that $\mathbf x_0\in \partial\Omega $ such that $ 	\Cor^{r_1}_{\mathbf x_0} \Subset \Omega $ or $ 	{\mathcal C}^{r_1}_{\mathbf x_0} \Subset \Omega $ where $r_1\in \mathbb R_+$, then there exists $r_0\in (0,r_1)$ such that  $ (\mathbf E^t,\mathbf H^t) \in C^{1/2}( \overline{ B_{r_0} (\mathbf x_0  )})^3 $. 
\end{lem}


\begin{proof}
Let $B_R$ be an open ball centered at the origin with the radius $R\in \mathbb R_+$ such that $\Omega\Subset B_R $. We first show that ${\mathbf E}^t\in H^1(B_R;\mathbb C^3 ) $ and ${\mathbf H}^t\in H^1(B_R;\mathbb C^3 ) $ by using Lemma \ref{lem:31 regu}. Considering  the electromagnetic  medium scattering problem \eqref{eq:Maxwell2},  it is ready to know that 
\begin{equation}\label{eq:trans B_R}
	\nu \wedge {\mathbf E}^t_-=\nu \wedge {\mathbf E}^t_+\quad \mbox{ on } \quad \partial B_R. 
\end{equation}
Since ${\mathbf E}^t_+$ is real analytic in $\mathbb R^3 \backslash B_R$, let $\mathbf G$ be the unique analytic continuation of ${\mathbf E}^t_+$ in $B_R$. By \eqref{eq:trans B_R}  one has
\begin{align}\label{eq:trans B_R G}
\nu \wedge {\mathbf E}^t_-=\nu \wedge {\mathbf G} \quad \mbox{ on } \quad \partial B_R \quad {and} \quad \mathbf G \in H^1(B_R; \mathbb C^3 ). 
\end{align}
By virtue of \eqref{eq:Maxwell2} and \eqref{eq:trans B_R G}, it arrives at
\begin{equation}\label{eq:regu lem}
	\left\{
	\begin{split}
	& \nabla\wedge  {\EM{E}}^t-\im\omega \Mp_0{\EM{H}}^t=\EM{J}_{m} \quad \mbox{in $B_R$},\\
	& \, \nabla\wedge {\EM{H}}^t+\im\omega\Ep_0 {\EM{E}}^t=\EM{J}_{e} \hspace{0.5cm} \mbox{in $B_R$},\\
	&\, \nor\cros {\EM{E}}^t=\nor\cros {\EM{G}} \hspace*{1.7cm}\mbox{on $\partial  B_R $, } 
	\end{split}\right.
	\end{equation}
where 
\begin{align}\notag
	\mathbf J_m=\begin{cases}
		( \varepsilon_0-\gamma(\mathbf x) ) \mathbf E^t &\mbox{in} \quad \Omega,\\
		\mathbf 0 &\mbox{in}\quad B_R\backslash \Omega,
	\end{cases}
	\quad \mathbf J_e=\begin{cases}( \mu(\mathbf x) -\mu_0) \mathbf H^t&\mbox{in} \quad \Omega,\\
	\mathbf 0 &\mbox{in}\quad B_R\backslash \Omega,
\end{cases}	
\end{align}
 and $\varepsilon_0,\mu_0$ are two positive constants. Due to \eqref{eq:adm para} and \eqref{eq:EH curl cond}, we immediately know that  $\mathbf J_e, \mathbf J_m\in L^2(B_R;\mathbb C^3 )$. By directly calculations, using \eqref{eq:adm para} and the divergence free property of $\gamma(\mathbf x)\mathbf E^t (\mathbf x)$ in $\mathbb R^3$, it  yields that
\begin{align}\notag
	\nabla \cdot \mathbf J_m=\varepsilon_0 \nabla \cdot \mathbf E^t=-\frac{\varepsilon_0 }{\gamma(\mathbf x) }\  \nabla \mathbf E^t\cdot \nabla \gamma(\mathbf x) \in L^2(\Omega;\mathbb C^3 ). 
	\end{align}
By virtue of Lemma \ref{lem:31 regu}	, one has ${\mathbf H}^t\in H^1(B_R;\mathbb C^3 ) $. Using similar arguments,  we know that ${\mathbf E}^t\in H^1(B_R;\mathbb C^3 ) $.

	Since  $ (\Omega ; \varepsilon ,\mu,\sigma ) $ is admissible,  we know that $\mu=\mu_1$ and $\gamma=\gamma_1$ are two positive constant in $B_{r_1}(\mathbf x_0) \cap \Omega$. Consider \eqref{eq:Maxwell2}, by direct calculations, we can deduce that
	\begin{equation}\nonumber
	\begin{cases}
		\nabla \wedge (\nabla \wedge \EM{E}^t)-\mathrm{i} \omega \mu_1 \EM{H}^t=\mathbf 0,\ \nabla \cdot(\nabla \wedge \EM{H}^t)+\mathrm{i} \omega\gamma_1(\nabla \cdot \EM{E}^t)=\mathbf 0\ &\mbox{ in } B_{r_1}(\mathbf x_0 )  \cap  \Omega,\\
				\nabla \wedge (\nabla \wedge \EM{E}^t)-\mathrm{i} \omega \mu_0 \EM{H}^t=\mathbf 0,\ \nabla \cdot(\nabla \wedge \EM{H}^t)+\mathrm{i} \omega\varepsilon_0(\nabla \cdot \EM{E}^t)=\mathbf 0\  &\mbox{ in } B_{r_1}(\mathbf x_0 )  \backslash \Omega. 
	\end{cases}
	\end{equation}
By  virtue of $ \nabla \cdot (\nabla \wedge \EM{H}^t)=\mathbf 0 $, it holds that
	\begin{equation}\label{eq:3}
\Delta \EM{E}^t+ k^2q\EM{E}^t=\mathbf 0 \quad \mbox{in\ \ $B_{r_1}(\mathbf x_0 )$},
	\end{equation}
	where $ k=\omega\sqrt{\mu_0\varepsilon_0} $ and
	$$ q=\chi_{B_{r_1}(\mathbf x_0 )  \cap \Omega   } \, \times \frac{\mu_1\gamma_1}{\mu_0\varepsilon_0} + \chi_{ B_{r_1}(\mathbf x_0 )  \backslash \Omega    } . 
	$$
	Due to ${\mathbf E}^t\in H^1(B_R;\mathbb C^3 ) $, by elliptic interior regularity and Sobolev embedding property, we can obtain that
	$ \mathbf E^t \in C^{\frac{1}{2}}(B_{r_0}(\mathbf x_0 )   )^{3} $, where $r_0\in (0, r_1)$. Similarly we can show that $ \mathbf H^t \in C^{\frac{1}{2}}(B_{r_0}(\mathbf x_0 )  )^{3} $.
	\end{proof}

\begin{thm}\label{thm:main7}
	Let $(\Omega ; \varepsilon,\mu,\sigma )$ and $(\Omega ; \varepsilon',\mu',\sigma' )$ be two admissible medium scatterers associated  with the electromagnetic medium scattering  \eqref{eq:Maxwell2}, where  $(\mathbf{E}^t, \mathbf{H}^t)$ and  $ ({\mathbf{E}^t}', {\mathbf{H}^t}')$ are the  corresponding total wave field to \eqref{eq:Maxwell2}. Suppose that  $\mathbf{E}_\infty$ and $\mathbf{E}_\infty'$ are the associated far fields. 
	 If \begin{equation}\label{eq:EE'}
		\mathbf{E}_\infty(\hat{\mathbf x} )=\mathbf{E}_\infty'(\hat{\mathbf x}) \quad \text{for all}\quad  \hat{\mathbf x}\in\mathbb{S}^2\ \text{and a fixed incident wave}\ (\mathbf E^i,\mathbf H^i),
	\end{equation}
	then the set difference
	$\Omega\Delta\Omega'$  defined in \eqref{eq:diff01} cannot possess a conical  corner. 

\end{thm}
\begin{proof}
	We prove this theorem by contradiction. Assume that there exists a conical corner $ \Cor _{\mathbf x_0}^{r_0} $ defined by \eqref{eq:coner0} such that 
	\begin{equation}\label{eq:coner}
		\Cor_{\mathbf x_0}^{r_0}\subset \overline \Omega \backslash \Omega'.
	\end{equation}
Suppose that $ (\mathbf E^t,\mathbf H^t) $ and $ ({{\mathbf E}^t}',{{\mathbf H}^t }^{'}) $ are the total wave field associated with the scatterers $ (\Omega ; \varepsilon,\mu,\sigma ) $ and $(\Omega ; \varepsilon',\mu',\sigma ' ) $ respectively. Hence, from \eqref{eq:coner}, we can have 
	\begin{equation}\label{eq:unique conical}
		\begin{cases}
		\nabla\wedge \EM{E}^t-\mathrm{i}\omega\Mp \EM{H}^t=\mathbf 0,\quad
		\nabla \wedge \EM{H}^t+\mathrm{i}\omega\Eci \EM{E}^t=\mathbf 0
		&\mbox{in\ \ $\Cor_{\mathbf x_0}^{r_0}$},\\
		\nabla\wedge {\EM{E}^t}^{'}-\mathrm{i}\omega\Mp_0 {\EM{H}^t}^{'}=\mathbf 0,\quad
		\nabla\wedge {\EM{H}^t}^{'}+\mathrm{i}\omega\Ep_0 \EM{E}_t^{'}=\mathbf 0
		&\mbox{in\ \ $\Cor_{\mathbf x_0}^{r_0}$},\\
		\end{cases}
	\end{equation}
	where $\gamma=\varepsilon+\mathrm{i}\sigma/\omega$. In view of \eqref{eq:EE'}, by Rellich theorem, we have $ (\mathbf E^t,\mathbf H^t)=({\mathbf E^t}^{'},{\EM{H}^t}^{'}) $ in $ \mathbf G $, where $ \mathbf G $ is defined in \eqref{eq:G}. Therefore, it readily to know that
	\begin{equation}\label{eq:boudary}
	\nu \wedge \mathbf E^t=\nu \wedge {\mathbf E^t}^{'},\ \nu \wedge \mathbf H^t=\nu \wedge {\mathbf H^t}^{'}\
			\mbox{on\ \ $\partial \Cor_{\mathbf x_0}^{r_0} \backslash \partial B_{r_0}(\mathbf x_0)$},
	\end{equation}
	where $\nu$ is the exterior unit normal vector to $  \partial \Cor_{\mathbf x_0,\theta_{0}}^{r_0} \backslash \partial B_{r_0}(\mathbf x_0) \subset \partial \mathbf G $.  Using the similar argument in the proof of Theorem \ref{thm:main5}, by virtue of \eqref{eq:unique conical} and \eqref{eq:boudary}, it is readily to see that  $(\widehat{\mathbf{E}},\widehat{\mathbf{H}}):=(\mathbf{E}^t,\mathbf{H}^t)-(\mathbf{E}',\mathbf{H}')$ satisfies 
	\begin{equation}\label{eq:rega3 unique}
	\left\{
	\begin{split}
	& \nabla\wedge  \widehat{\EM{E}}-\im\omega \Mp_0 \widehat{\EM{H}}=\EM{J}_{1} \quad \mbox{in $\Cor^{r_0}_{\mathbf x_0}$},\\
	& \, \nabla\wedge \widehat{\EM{H}}+\im\omega\Ep_0 \widehat{\EM{E}}=\EM{J}_{2} \quad \mbox{in $\Cor^{r_0}_{\mathbf x_0}$},\\
	&\, \nor\cros \widehat{\EM{E}}=\nor\cros\widehat{\EM{H}}=\mathbf 0 \hspace*{0.8cm}\mbox{on $\partial \Cor^{r_0}_{\mathbf x_0}\setminus \partial B_{r_0}(\mathbf x_0)$ ,} 
	\end{split}\right.
	\end{equation}
	where
	$\mathbf{J}_1= \mathrm i \omega (\mu(\mathbf x)-\mu_0)\mathbf{H}^t\quad\mbox{and}\quad\mathbf{J}_2=\mathrm i \omega (\gamma(\mathbf x)-\varepsilon_0)\mathbf{E}^t. $


According to Lemma \ref{lem:36 max}, we know that
	$ \mathbf E^t \in C^{\frac{1}{2}}( \overline{ \Cor_{ \mathbf x_0}^{r_0} })^{3} $. Similarly we can show that $ \mathbf H^t \in C^{\frac{1}{2}}(\overline{ \Cor_{\mathbf x_0}^{r_0} } )^{3} $.Therefore, in view of \eqref{eq:rega3 unique}, by Theorem~\ref{thm:main0}, we have $ \mathbf E^t(\mathbf x_0)=\mathbf H^t(\mathbf x_0)=\mathbf 0 $, which contradicts to the condition (3) in Definition~\ref{def:adm}.

	The proof is complete.
\end{proof}

We proceed to prove a global unique identifiability result for the shape determination of an admissible scatterer of coronal shape described by Definition~\ref{def:finite c} under certain priori knowledge on the underlying scatterer.

\begin{thm}\label{thm:34}
	Suppose that $\Omega=\cup _{j=1}^{\ell} \Cor_{\mathbf x_j,\theta_{j}} \cup D $ and $\Omega'=\cup _{j=1}^{\ell'} \Cor_{\mathbf x_j^{'},\theta_{j}^{'}} \cup D' $ be two scatters of coronal shape  described by Definition \ref{def:finite c}. Let $(\Omega ; \varepsilon,\mu,\sigma )$  and $(\Omega ; \varepsilon',\mu',\sigma ')$  be two admissible electromagnetic medium  scatters associated with the electromagnetic medium scattering  \eqref{eq:Maxwell2} and the incident wave $ (\mathbf E^i,\mathbf H^i) $. If \eqref{eq:EE'} and 
	\begin{equation}\label{eq:thm34 cond}
		D=D',\quad \theta_{j}=\theta_{i}^{'}\quad   \mbox{for} \quad  j\in \{1,\ldots,\ell\}\quad \mbox{and}  \quad  i\in \{1,\ldots,\ell'
		\}\quad \mbox{fulfilling} \quad \mathbf x_j=\mathbf x_i^{'},
	\end{equation}
	are satisfied, then $\Omega=\Omega'$, 
	 $\varepsilon (\mathbf  x_j)=\varepsilon '(\mathbf  x_j),  \sigma(\mathbf  x_j)=\sigma'(\mathbf  x_j)$, and $ \mu(\mathbf  x_j)=\mu'(\mathbf  x_j)$
 for $j\in {1,...,\ell}$.
 
 If  \eqref{eq:EE'}, $	D=D'$ and 
	\begin{equation}\notag 
	 \varepsilon ({\mathbf x_j})\not=\varepsilon( \mathbf x_{i}^{'})   \mbox{ or }  \gamma({\mathbf x_j})\not=\gamma( \mathbf x_{i}^{'})    \mbox{ for }  j\in \{1,\ldots,\ell\} \mbox{ and }    i\in \{1,\ldots,\ell'
		\} \mbox{ fulfilling }  \mathbf x_j=\mathbf x_i^{'},
	\end{equation}
	are satisfied, where $\gamma$ is defined in \eqref{eq:Maxwell2}, then $\Omega=\Omega'$. 
\end{thm}
\begin{proof}
	We prove this theorem by contradiction. Suppose that $\ell\neq \ell'$; or if $\ell=\ell'$ but $\mathbf x_j\neq  \mathbf x_j'$ or $\theta_{j}=\theta'_{j}$ ($j=1,\ldots, \ell$), in view of \eqref{eq:thm34 cond}, without loss of generality, one can
claim that there must exist a conical corner $ \Cor_{\mathbf  x_{j_{0}}}^{r_{0}} $ defined by \eqref{eq:coner0} such that $ \Cor_{\mathbf  x_{j_{0}}}^{r_{0}} \in \bar{\Omega} \backslash \Omega'$. Since   \eqref{eq:EE'} is satisfied, by virtue ofTheorem \ref{thm:main5},  we  directly  get the contradiction.

Since $\Omega=\Omega'$, with the help of \eqref{eq:EE'} and Rellich theorem, it holds that
\begin{equation}\label{eq:gama mu 1}
		\begin{cases}
		\nabla\wedge \EM{E}^t-\mathrm{i}\omega\Mp \EM{H}^t=\mathbf 0,\quad
		\nabla \wedge \EM{H}^t+\mathrm{i}\omega\Eci \EM{E}^t=\mathbf 0
		&\mbox{in\ \ $\mathcal K_{\mathbf x_j}^{r_j}$},\\
		\nabla\wedge {\EM{E}^t}^{'}-\mathrm{i}\omega\Mp' {\EM{H}^t}^{'}=\mathbf  0,\quad
		\nabla\wedge {\EM{H}^t}^{'}+\mathrm{i}\omega\gamma' \EM{E}_t^{'}=\mathbf 0
		&\mbox{in\ \ $\mathcal K_{\mathbf x_j}^{r_j}$},\\
		\quad\nu \wedge \mathbf E^t=\nu \wedge {\mathbf E^t}^{'},\ \nu \wedge \mathbf H^t=\nu \wedge {\mathbf H^t}^{'}\
			&\mbox{on\ \ $\partial \mathcal K_{\mathbf x_j}^{r_j} \backslash \partial B_{r_j}(\mathbf x_j)$},
		\end{cases}
	\end{equation}		
	where $\mathbf x_j$ is the apex of  the conic corner $\mathcal K_{\mathbf x_j}^{r_0}$ 	and $\nu$ is the exterior unit normal vector to $\Omega$.  It can verify that \eqref{eq:gama mu 1}  can be written as
	\begin{equation}\label{eq:gama mu}
		\begin{cases}
		\nabla\wedge  \widetilde {\EM{E}} -\mathrm{i}\omega\Mp' \widetilde{ \EM{H}}=\widetilde{ \mathbf J}_1,\quad
		&\mbox{in\ \ $\mathcal K_{\mathbf x_j}^{r_j}$},\\
		\nabla \wedge \widetilde{ \EM{H}}+\mathrm{i}\omega\Eci ' \widetilde{ \EM{E} }=\widetilde {\mathbf J}_2 \quad 
		&\mbox{in\ \ $\mathcal K_{\mathbf x_j}^{r_j}$},\\
		\quad\nu \wedge\widetilde{  \mathbf E}=\nu \wedge \widetilde{ {\mathbf H}}=\mathbf 0,\quad
			&\mbox{on\ \ $\partial \mathcal K_{\mathbf x_j}^{r_j} \backslash \partial B_{r_j}(\mathbf x_j)$},
		\end{cases}
	\end{equation}		
	where  $\widetilde{\mathbf E}= {\EM{E}^t}- {\EM{E}^t}^{'}$, $ \widetilde{\mathbf H}= {\EM{H}^t}- {\EM{H}^t}^{'} $, $\widetilde{ \mathbf J}_1=\mathrm{i}\omega (\mu-\mu')\mathbf{H}^t $	 and $\widetilde{ \mathbf J}_2=\mathrm{i}\omega (\gamma'-\gamma)\mathbf{E}^t $. According to Lemma \ref{lem:36 max}, we know that  $\widetilde{ \mathbf J}_1(\mathbf x)\in C^{1/2} ( \overline{ \mathcal C_{\mathbf x_j}^{r_j}  } )^3  $ and $\widetilde{ \mathbf J}_2(\mathbf x) \in C^{1/2} ( \overline{ \mathcal C_{\mathbf x_j}^{r_j} } )^3   $.  By using the similar argument of Theorem \ref{thm:main0}, we can prove that $\widetilde{ \mathbf J}_1(\mathbf x_j)=\widetilde{ \mathbf J}_2(\mathbf x_j)=\mathbf 0$, which implies that $\gamma (\mathbf x_j)=\gamma'(\mathbf  x_j)  $ and $\mu(\mathbf  x_j)=\mu'(\mathbf  x_j)$ by noting $\mathbf{E}^t (\mathbf x_j)\neq \mathbf 0$ and $\mathbf{H}^t (\mathbf x_j)\neq \mathbf 0$.

	
	The second part of this theorem can be proved by using a similar argument  for proving $\Omega=\Omega'$ under the condition\eqref{eq:J1J2} in Theorem \ref{thm:global source}.

	The proof is complete.
\end{proof}

\begin{rem}
	Similar to Theorem  \ref{thm:34},  by virtue  of the regularity results on  the total wave field associated with \eqref{eq:Maxwell2}  near a  polyhedral corner in Lemma \ref{lem:36 max},  utilizing the local uniqueness results from   \cite[Theorem 4.3]{BLX2020} with respect to a polyhedral corner, one can establish a global unique determination for the shape by a single far field measurements for an admissible convex polyhedron medium scatterer $(\Omega ; \varepsilon,\mu,\sigma )$.  Using a similar argument inn Theorem \ref{thm:34}, we can also prove  $\varepsilon (\mathbf  x_j)$, $\sigma(\mathbf  x_j)$ and  $  \mu(\mathbf  x_j)$ ($\forall \mathbf  x_j \in \mathcal V(\Omega ) 
	$),  where $ \mathcal V(\Omega )$ is the set of the vertexes of $\Omega$,   can  be uniquely determined by a single far field measurement.
\end{rem}

By virtue of the contradiction argument, Theorem \ref{thm:main5} and Lemma \ref{lem:36 max}, we have the following theorem  which indicates that an electromagnetic medium possessing a conical corner always scatters. 
\begin{thm}\label{thm:medium conical scat}
	Consider the electromagnetic medium scattering problems \eqref{eq:Maxwell2}. Let $(\Omega ; \varepsilon,\mu,\sigma )$   be the medium scatterer  associated with \eqref{eq:Maxwell2},  where  $(\mathbf{E}^t, \mathbf{H}^t)$ is the  corresponding total wave field to \eqref{eq:Maxwell2}. Suppose that the condition \eqref{eq:adm para} for the physical parameter $\varepsilon$, $\mu$ and $\sigma$ is fulfilled. If $\Omega$ has a a conical corner $\mathcal K_{\mathbf x_0}^{r_0}$ described by \eqref{eq:coner0} and  the total wave field $\mathbf E^t$ or $\mathbf H^t$ is non-vanishing at $\mathbf x_0$ in the sense of \eqref{eq:nn},  then $\Omega$ always scatters for any incident wave satisfying \eqref{eq:Maxwellh}. 
	\end{thm}

\section{Spectral geometry of transmission eigenfunctions}\label{sec:4}



In this section, we consider the following transmission eigenvalue problem: 
\begin{equation}\label{eq:ITP}
		\begin{cases}
			\nabla\wedge \EM{E}^t-\mathrm{i}\omega\Mp \EM{H}^t=\mathbf 0,\quad
			\nabla \wedge \EM{H}^t+\mathrm{i}\omega\Eci \EM{E}^t=\mathbf 0
			&\mbox{in\ \ $\Omega$},\\
			\nabla\wedge \EM{E}^0-\mathrm{i}\omega\Mp_0 \EM{H}^0=\mathbf 0,\quad
			\nabla\wedge \EM{H}^0+\mathrm{i}\omega\Ep_0 \EM{E}^0=\mathbf 0
			&\mbox{in\ \ $\Omega$},\\
			\quad\nu\cros \EM{E}^t=\nu\cros \EM{E}^0,\quad
			\qquad\nu\cros \EM{H}^t=\nu\cros \EM{H}^0
			&\mbox{on\ \ $\Gamma $},
		\end{cases}
	\end{equation}
where $\Gamma\Subset \partial \Omega$ and $(\mathbf{E}^t,\mathbf{H}^t)\in H(\mathrm{curl}, \Omega)\times H(\mathrm{curl}, \Omega )$ and $(\mathbf{E}^0,\mathbf{H}^0)\in H(\mathrm{curl}, \Omega)\times H(\mathrm{curl}, \Omega)$. If there exists $\omega\in\mathbb{R}_+$ such that there exist nontrivial solutions $(\mathbf{E}^t,\mathbf{H}^t)\in H(\mathrm{curl}, \Omega)\times H(\mathrm{curl}, \Omega )$ and $(\mathbf{E}^0,\mathbf{H}^0)\in H(\mathrm{curl}, \Omega)\times H(\mathrm{curl}, \Omega)$ satisfying \eqref{eq:ITP}, then $\omega$ is called an interior transmission eigenvalue and $(\mathbf{E}^t, \mathbf{H}^t), (\mathbf{E}^0, \mathbf{H}^0)$ are called the corresponding transmission eigenfunctions (cf. \cite{CHreview,Liu22}). According to \cite{Liu22}, if $\Gamma=\partial\Omega$, \eqref{eq:ITP} is referred to as the full-data transmission eigenvalue problem, otherwise it is referred to as the partial-data transmission eigenvalue problem. It is clear by \eqref{eq:ITP old} that if invisibility occurs for the scattering problem \eqref{eq:Maxwell2}, then the total and incident fields form a pair of (full-data) transmission eigenfunctions.

In view of the proof of Theorem \ref{thm:main5}, one readily has  the following theorem on the locally vanishing characterization of transmission eigenfunctions to \eqref{eq:ITP} near a conical corner when they have H\"older continuity near the concerned conical corner.     
\begin{thm}\label{thm:main trans}
	Let $(\mathbf{E}^t, \mathbf{H}^t)$ and $(\mathbf{E}^0, \mathbf{H}^0)$ be a pair of eigenfunctions to the interior transmission eigenvalue function \eqref{eq:ITP} associated with the eigenvalue $\omega\in\mathbb{R}_+$. 
	Assume that $\Omega$ possesses a conical corner $\mathcal{K}_{r_0}$ with $\mathbf 0\in \Gamma$, where $ \mathbf 0 $ is the apex of the cone $ \mathcal{K} $, and 
	\begin{equation}\label{eq:rega1}
	(\mu-\mu_0)\mathbf{H}^t,\ (\gamma-\varepsilon_0)\mathbf{E}^t\in C^\alpha(\overline{\mathcal{K}_{r_0}})^3,
	\end{equation} 
	for some $\alpha\in (0, 1)$. Then there holds
	\begin{equation}\label{eq:rega2}
	(\mu(\mathbf 0)-\mu_0)\mathbf{H}^t(\mathbf 0)=(\gamma(\mathbf 0)-\varepsilon_0)\mathbf{E}^t(\mathbf 0)=\mathbf 0. 
	\end{equation}
	Furthermore, if $\mu(\mathbf 0)\neq \mu_0\neq 0$ and $\gamma(\mathbf 0)\neq \varepsilon_0$, one has $\mathbf{H}^t(\mathbf 0)=\mathbf{E}^t(\mathbf 0)=\mathbf 0.$
\end{thm}

By our study in Section~\ref{sec:3}, we see that the H\"older continuity is a physically unobjectionable condition and can be satisfied for certain practical scenarios. This is crucial in establishing the unique recovery result in Subsection~\ref{subsec:in me}. We proceed to consider the vanishing property of the transmission eigenfunctions when the H\"older continuity condition is weakened. In fact, it is shown in \cite{2021} through numerics that the vanishing property may hold under more general regularity conditions on the underlying transmission eigenfunctions. In the rest of this section, we shall show that under a regularity condition in terms of the Herglotz approximation, the locally vanishing property still holds for electromagnetic transmission eigenfunctions. In the context of acoustic transmission eigenvalue problem, it is sharply quantified in \cite{LT1} that the regularity condition in terms of the Herglotz approximation is indeed weaker than the H\"older continuity. 

In what follows, we let $ \mu $ and $ \gamma  $ in \eqref{eq:ITP} be positive  constants. Next, we introduce some auxiliary results concerning the electromagnetic Herglotz approximation. 

Let an electric Herglotz wave function with the kernel $ \mathbf g$ be defined by 
	\begin{equation}\label{eq:Eg}
	\mathbf{E}_{\mathbf g}(\mathbf x)=\int_{ \mathbb S^{2}}\mathbf g(\mathbf d)\exp(\im k_1 \mathbf x\cdot \mathbf d) \quad \text{for all}\  \mathbf x \in \mathbb R^3, \quad k_1=\omega\sqrt{\mu \gamma}, 
	\end{equation}
	where $ k_{1} \in \mathbb R_{+} $, $\mathbf d\in \mathbb S^2$ and $ \mathbf g\in L^2(\mathbb  S^{2}) $. Hence one can easily verify that  $\mathbf{E}_{\mathbf g}$  is an entire solution to 
	\begin{equation}\notag
	\nabla \wedge \nabla \wedge \mathbf{E}-k_{1}^2\mathbf{E}=\mathbf 0,
	\end{equation}
	where $ k_{1}$ is defined in \eqref{eq:Eg}.  Similarly, we define 
	\begin{equation}\label{eq:Hf}
	\mathbf{H}_{\mathbf f}= \frac{1}{\im \omega \mu_0}\nabla \wedge \mathbf{E}_{\mathbf g},
	\end{equation}
	where $ \mathbf{E}_{\mathbf g} $ is given by \eqref{eq:Eg}. Therefore one has
	$$
	\mathbf{H}_{\mathbf f}=\int_{ \mathbb S^{2}}\mathbf f(\mathbf d)\exp(\im k_1 \mathbf x\cdot \mathbf d) \mathrm {d} \sigma({\mathbf d})\quad \text{for all}\  \mathbf x \in \mathbb R^3,
	$$
	where $\mathbf f(\mathbf d)=\frac{k_1}{\omega \mu_0}\mathbf d\wedge \mathbf g(\mathbf d)$ is the kernel of $\mathbf{H}_{\mathbf f}$. $\mathbf{H}_{\mathbf f}$ is said to be 
	a magnetic Herglotz wave function with the kernel  $\mathbf f$.  Similarly,    $\mathbf{H}_{\mathbf f}$ is an entire solution to 	$\nabla \wedge \nabla \wedge \mathbf{H}-k_{1}^2\mathbf{E}=\mathbf 0$.

\begin{lem}\label{lem:approx}\cite{monk}
Let $ D $ be a bounded Lipschitz domain with a connected  
	complement. Then the set of electric Herglotz wave functions $ \mathbf{E}_{\mathbf g} $ with the form \eqref{eq:Eg}   is dense with respect to the $ H({\rm curl};D) $ norm in the set of solutions to 
	\begin{equation}
	\nabla \wedge \nabla \wedge \mathbf E^t-k_{1}^2\mathbf E^t= \mathbf 0
	\end{equation}
	in $ D $, where $ k_{1}$ is defined in \eqref{eq:Eg}. Similarly, let the magnetic  Herglotz wave function $ \mathbf{H}_{\mathbf f} $ be defined by   \eqref{eq:Hf}, where the electric Herglotz wave functions $ \mathbf{E}_{\mathbf g} $ is given by \eqref{eq:Eg}. The magnetic  Herglotz wave function $ \mathbf{H}_{\mathbf f} $ can approximate any solution $\mathbf H \in H(\mathrm{curl},D)$ satisfying $\nabla \wedge \nabla \wedge \mathbf{E}-k_{1}^2\mathbf{E}=\mathbf 0$ in $D$ (in the distribution sense) with arbitrary accuracy,
\end{lem}

The following lemma can be summarized from the the proof of Theorem 1.1 of \cite{BLX2020}.
\begin{lem}\label{lem:wv}
Recall that $ \EM{V}(\mathbf x) $ and $ \EM{W}(\mathbf x) $ are defined in \eqref{eq:vw}, $ \nu \in \mathbb S^2 $, $\pare{\EM{E},\EM{H}}\in H(\mathrm{curl},\Omega)\times H(\mathrm{curl},\Omega)$ is a solution to the Maxwell system \eqref{eq:MaxwellO}, and $ \Cor_{r_{0}} $ is the truncated conical  cone  defined in \eqref{eq:coner0}. Assume that $\mathbf 0\in \partial \Omega$, and  there exists $r_0\in \mathbb R_+$ such that ${\mathcal K}_{r_0} \Subset \Omega \cap {\mathcal K}$. It holds that
	
	\begin{equation}\label{eq:lem26}
	\begin{split}
	\left|\int_{\partial \Cor_{r_{0}}\cap \partial B_{r_0}}
	\EM{W}\cdot\pare{\nor\cros\EM{E}} \mathrm{d} \mathbf {x}\right| &\le 3C_{\Cor,r_0}k^2e^{-\delta r_0\tau}
	\|\EM{E}\|_{H(\curl,\, \Cor_{r_{0}})},\\
	\left|\int_{\partial \Cor_{r_{0}}\cap \partial B_{r_0}}
	\EM{V}\cdot\pare{\nor\cros\EM{H}} \mathrm{d} \mathbf {x}\right| &\le 3C_{\Cor,r_0}k^2e^{-\delta r_0\tau}
	\|\EM{H}\|_{H(\curl,\, \Cor_{r_{0}})}
	\end{split}
	\end{equation}
	as $\tau \rightarrow +\infty$,  where  $ C_{\Cor,r_0} $ is a positive constant only depending on the surface measure of   $\partial \mathcal K_{r_0} \cap \partial B_{r_0}$. 
\end{lem}

Several auxiliary lemmas in the following can be obtained from the proof of \cite[Lemmas 2.4-2.7]{2021}, which play important roles in deriving our main result in Theorem~\ref{th:hcurl} in what follows.

\begin{lem}\label{lem:v}
	Recall that $ \Cor_{r_0} $ and $ \EM{V}(\mathbf x) $ are defined in \eqref{eq:coner0} and \eqref{eq:vw} respectively, where $\boldsymbol{\rho}$ and $\mathbf p$ in $\mathbf V(\mathbf x)$ satisfy \eqref{eq:r} and 
	\eqref{eq:vw}. 
	Then there holds that 
	\begin{equation}\label{eq:vkr0}
	\left| \int_{\Cor_{r_{0}}} \EM{V}(\mathbf x) {\mathrm d} \mathbf x\right|\le C\tau ^{-3}(1+\mathcal{O}(\tau^{-2} )),
	\end{equation}
	as $\tau \rightarrow +\infty$, where $ C $ is a positive constant only depending on $\theta_0$ and $c_{\Cor} $, $c_{\Cor}$ is given by \eqref{eq:innerdot}.
\end{lem}

%
%

\begin{lem}\label{lem:vl2}
	Recall that $ \Cor_{r_0} $ and $ \EM{V}(\mathbf x) $ are defined in \eqref{eq:coner0} and \eqref{eq:vw} respectively,  where $\boldsymbol{\rho}$ and $\mathbf p$ in $\mathbf V(\mathbf x)$ satisfy \eqref{eq:r}. Then there holds that
	\begin{equation}\label{eq:vl2}
	\left\|\EM{V}(\mathbf x) \right\|_{L^2(\Cor_{r_0})} \le C \tau^{-\frac{3}{2}}(1+\mathcal{O}(\tau^{-2} )), 
	\end{equation}
	as $\tau \rightarrow +\infty$, where $ C $ is a positive constant only depending on $\theta_0$ and $\delta $, $\delta$ is given by \eqref{eq:innerdot}.
\end{lem}

\begin{lem}\label{lem:egv}
	Assume that  $ \EM{V} $ is defined in \eqref{eq:vw} . Let $ \mathbf{E}_{\mathbf  g_j} $ be defined by
	\begin{equation}\label{eq:Egj}
	\mathbf{E}_{\mathbf g_j}(\mathbf x)=\int_{ \mathbb S^{2}}\mathbf g_j(\mathbf d)\exp(\im k_1 \mathbf x\cdot \mathbf d),
	\end{equation}
	where $\mathbf g_j\in L^2(\mathbb S^2)$ ,  $ k_{1}=\omega\sqrt{\varepsilon\gamma} $ and the constants $ \varepsilon, \gamma \in \mathbb R_{+} $.  Then $\mathbf{E}_{\mathbf g_j} \in C^1 (\mathcal K_{r_0} )^3$,  $r_0\in \mathbb R_+$ and $\mathcal K_{r_0}$  is defined by \eqref{eq:coner0} and has the expansion
	\begin{equation}\label{eq:Ejex}
	\mathbf{E}_{\mathbf g_j}(\mathbf x)=\mathbf{E}_{\mathbf g_j}(\mathbf 0)+ \delta \mathbf{E}_{\mathbf  g_{j}}(\mathbf x),\quad | \delta \mathbf{E}_{\mathbf  g_{j}} (\mathbf x)|\leq  \|  \mathbf{E}_{\mathbf  g_{j}}\|_{ C^1 (\mathcal K_{r_0} )^3 } |\mathbf x|. 
	\end{equation}
	Suppose that $ \mathbf g_{j} $ satisfies the following condition
	\begin{equation}\label{eq:g}
	\ \left\| \mathbf g_{j} \right\| _{L^2(\Cor_{r_{0}})^{3}}\le j^{\beta},
	\end{equation} 
	where $ \beta \in  \mathbb R_{+} $,  we have
	\begin{equation}\label{eq:lem35}
	\left| \int_{\Cor_{r_{0}}} (\gamma-\varepsilon_0)\delta \mathbf{E}_{\mathbf  g_{j}}\cdot\EM{V}  \right| \le Cj^{\beta}\tau^{-4}(1+ \mathcal O (  \tau^{-2} ))
	\end{equation}
	as $\tau \rightarrow +\infty$,  where $ C $ is a positive constant only depending on $\theta_0$ and $\delta $, $\delta $ is given by \eqref{eq:innerdot}. Similarly, let
	\begin{equation}\label{eq:35 Hj}
	\mathbf{H}_{\mathbf f_j}(\mathbf x)= \int_{ \mathbb S^{2}}\mathbf f_j(\mathbf d)\exp(\im k_1 \mathbf x\cdot \mathbf d) \mathrm {d} \sigma({\mathbf d}),\quad \mathbf f_j(\mathbf d)=\frac{k_1}{\omega \mu_0}\mathbf d\times \mathbf g_j(\mathbf d), 
	\end{equation}
	where $\mathbf g_j(\mathbf d)$ is given by \eqref{eq:Egj}.
	 It holds that $\mathbf{H}_{\mathbf f_j}(\mathbf x)=\frac{1}{\im \omega \mu_0}\nabla \times \mathbf{E}_{\mathbf g_j}$, where $\mathbf{E}_{\mathbf g_j}$ is defined in \eqref{eq:Egj}. Furthermore, we have $\mathbf{H}_{\mathbf f_j} \in C^1 (\mathcal K_{r_0} )^3$ and the expansion
	\begin{equation}\label{eq:Hjex}
	\mathbf{H}_{\mathbf f_j}(\mathbf x)=\mathbf{H}_{\mathbf f_j}(\mathbf 0)+ \delta \mathbf{H}_{\mathbf  f_{j}}(\mathbf x),\quad | \delta \mathbf{H}_{\mathbf  f_{j}} (\mathbf x)|\leq  \|  \mathbf{H}_{\mathbf  f_{j}}\|_{ C^1 (\mathcal K_{r_0} )^3 } |\mathbf x|. 
	\end{equation} 
	Suppose that $ \mathbf f_{j} $ satisfies the following condition
	\begin{equation}\label{eq:f}
	\ \left\| \mathbf f_{j} \right\| _{L^2(\Cor_{r_{0}})^{3}}\le j^{\beta '},
	\end{equation} 
	where $ \beta ' \in  \mathbb R_{+} $,  then we can obtain that
	\begin{equation}\label{eq:lem35 H}
	\left| \int_{\Cor_{r_{0}}}\delta \mathbf{H}_{\mathbf  f_{j}}\cdot\EM{V} \mathrm d \mathbf x \right| \le Cj^{\beta '}\tau^{-4}(1+ \mathcal O (  \tau^{-2} ))
	\end{equation}
	as $\tau \rightarrow +\infty$,  where $ C $ is a positive constant only depending on $\theta_0$ and $\delta $, $\delta$ is given by \eqref{eq:innerdot}. 
\end{lem}

The vanishing characterization of electromagnetic transmission eigenfunctions near a conical corner under a certain Herglotz wave approximation property can be obtain in the following theorem. 
\begin{thm}\label{th:hcurl}
	Suppose that $(\mathbf{E}^t, \mathbf{H}^t)$ and $(\mathbf{E}^0, \mathbf{H}^0)$ are a pair of eigenfunctions to the interior transmission eigenvalue  \eqref{eq:ITP} associated with the eigenvalue $\omega\in\mathbb{R}_+$. Assume $ \mu $ and $ \gamma  $ in \eqref{eq:ITP} are positive  constants,   $\mathbf 0 \in \Gamma$ and $ \Omega $ possesses a conical corner such that $\Omega \cap B_{r_0}(\mathbf 0)=\mathcal{K}_{r_0}$. Moreover, let   $ \mathbf{E}_{\mathbf g_{j}} $  and $ \mathbf{H}_{\mathbf f_{j}} $  be defined by \eqref{eq:Egj} and \eqref{eq:35 Hj} respectively, where the transmission eigenfunctions $\mathbf{E}^t, \ \mathbf{H}^t $ can be approximated by the electric and magnetic  Herglotz wave function   $ \mathbf{E}_{\mathbf g_{j}} $ and  $ \mathbf{H}_{\mathbf f_{j}} $ in the $ H(\curl;\ \Cor_{r_{0}})  $ norm, respectively with the approximation property 
	\begin{equation}\label{eq:eeg}
	\left\|\mathbf{E}^t-\mathbf{E}_{\mathbf g_{j}} \right\|_{H(\curl;\ \Cor_{r_{0}})^3} \le j^{-\zeta},\ \left\| \mathbf g_{j} \right\| _{L^2(\Cor_{r_{0}})^3}\le j^{\beta},
	\end{equation} 
	\begin{equation}\label{eq:heg}
	\left\|\mathbf{H^t}-\mathbf{H}_{\mathbf f_{j}} \right\|_{H(\curl;\ \Cor_{r_{0}})^3} \le j^{-\zeta'},\ \left\| \mathbf f_{j} \right\| _{L^2(\Cor_{r_{0}})^3}\le j^{\beta'},
	\end{equation}
		where $ \zeta ,\beta, \zeta', \beta'   \in  \mathbb R_{+} $ are fixed with $ \beta <\frac{2}{3} \zeta $ and  $ \beta' <\frac{2}{3} \zeta' $.
	It holds that 
	\begin{subequations}
		\begin{align}
\gamma \neq \varepsilon_0 \quad \mbox{implies}	\quad	&\lim_{\rho \to +0}\frac{1}{m(B(\mathbf 0,\rho) \cap \Cor_{r_{0}})}\int_{ B(\mathbf 0,\rho) \cap \Cor_{r_{0}}} | \mathbf E^t(\mathbf x) | \mathrm d \mathbf x=0, \label{eq:thm1} \\
\mu \neq \mu_0 \quad  \mbox{implies}\quad		&\lim_{\rho \to +0}\frac{1}{m(B(\mathbf 0,\rho) \cap \Cor_{r_{0}})}\int_{ B(\mathbf 0,\rho) \cap \Cor_{r_{0}}} |\mathbf H^t(\mathbf x) | \mathrm d \mathbf x=0.  \label{eq:thm2}
		\end{align}
	\end{subequations}
\end{thm}
\begin{proof}
It can be direct to see that  $ \EM{J}_{1} $ and $ \EM{J}_{2} $ in \eqref{eq:rega3} can be rewritten as 
	\begin{equation}\label{eq:new j1j2}
	\begin{split}
	\EM{J}_{1}&=(\mu-\mu_0)\mathbf{H}^t-(\mu-\mu_0)\mathbf{H}_{\mathbf g_{j}}+(\mu-\mu_0)\mathbf{H}_{\mathbf g_{j}}\\
	\EM{J}_{2}&=(\gamma-\varepsilon_0)\mathbf{E}^t-(\gamma-\varepsilon_0)\mathbf{E}_{\mathbf g_{j}}+(\gamma-\varepsilon_0)\mathbf{E}_{\mathbf g_{j}}. 
	\end{split}
	\end{equation}

	
	
	
	Since $\Cor_{r_{0}}$ is a convex conical corner with the apex $\mathbf 0$, we can choose $\mathbf d=(0,0,-1)^\top$ such that \eqref{eq:innerdot} is fulfilled, which implies that $\mathbf d^{\perp}=(\cos\varphi ,\sin \varphi , 0)^\top$ with $\varphi\in (0,2\pi]$.
 Therefore, for any $\pare{\EM{V},\EM{W}}$ defined by \eqref{eq:vw}, according to  Lemma~\ref{lem:Integral} and the boundary condition in \eqref{eq:rega3},  noting \eqref{eq:new j1j2}, there holds that 
	\begin{equation}\label{eq:herg1}
	\begin{split}
	\int_{\Cor_{r_{0}}} (\gamma-\varepsilon_0)(\mathbf{E}^t-\mathbf{E}_{\mathbf g_{j}})\cdot \EM{V}\mathrm d \mathbf x
	+\int_{\Cor_{r_{0}}} (\gamma-\varepsilon_0)\mathbf{E}_{\mathbf g_{j}}\cdot\EM{V}\mathrm d \mathbf x
	=&\int_{\partial \Cor_{r_{0}}\cap \partial B_{r_0}}
	\EM{W}\cdot\pare{\nor\cros\EM{E}}\mathrm d \sigma \\
	+&\int_{\partial \Cor_{r_{0}}\cap \partial B_{r_0}} \EM{V}\cdot\pare{\nor\cros\EM{H}}\mathrm d \sigma ,
	\end{split}
	\end{equation}	
Due to $ \mathbf{E}_{\mathbf g_{j}}  \in C^1(\overline{\mathcal{K}_{r_0}})^3 $, we substitute \eqref{eq:Ejex} into \eqref{eq:herg1}, it yields that
	\begin{equation}\label{eq:equality}
	\begin{split}
	(\gamma-\varepsilon_0)\int_{\Cor_{r_{0}}} \mathbf{E}_{\mathbf g_{j}}(0)\cdot \EM{V}\mathrm d \mathbf x
	=&\int_{\partial \Cor_{r_{0}}\cap \partial B_{r_0}}
	\EM{W}\cdot\pare{\nor\cros\EM{E}}\mathrm d \mathbf x+\int_{\partial \Cor_{r_{0}}\cap \partial B_{r_0}}\EM{V}\cdot\pare{\nor\cros\EM{H}}\mathrm d \mathbf x\\
	-&\int_{\Cor_{r_{0}}} (\gamma-\varepsilon_0)(\mathbf{E}^t-\mathbf{E}_{\mathbf g_{j}})\cdot\EM{V}\mathrm d \mathbf x-\int_{\Cor_{r_{0}}} (\gamma-\varepsilon_0)\delta \mathbf{E}_{\mathbf g_{j}}\cdot\EM{V}\mathrm d \sigma \\
	:=&I_{3}+I_{4}+I_{5}+I_{6},
	\end{split}
	\end{equation}
	where 
	\begin{equation}\notag
\begin{split}
		I_{3}&=\int_{\partial \Cor_{r_{0}}\cap \partial B_{r_0}}
		\EM{W}\cdot\pare{\nor\cros\EM{E}}\mathrm d \sigma ,\ I_{4}=\int_{\partial \Cor_{r_{0}}\cap \partial B_{r_0}}\EM{V}\cdot\pare{\nor\cros\EM{H}}\mathrm d \sigma ,\\
		I_{5}&=-\int_{\Cor_{r_{0}}} (\gamma-\varepsilon_0)(\mathbf{E}^t-\mathbf{E}_{\mathbf g_{j}})\cdot\EM{V}\mathrm d \mathbf x, \, 
		I_{6}=-\int_{\Cor_{r_{0}}} (\gamma-\varepsilon_0)\delta \mathbf{E}_{\mathbf g_{j}}\cdot\EM{V}\mathrm d \mathbf x.
\end{split}
	\end{equation}
	Using the Cauchy-Schwarz inequality, by virtue of \eqref{eq:vl2}  and \eqref{eq:eeg}, we can deduce that
	\begin{equation}\label{eq:331}
	\begin{split}
	|I_5|
		&\le
	\left| \gamma-\varepsilon_0\right| \left\| \mathbf{E}^t-\mathbf{E}_{\mathbf g_{j}} \right\|_{L^2} \cdot \left\|\EM{V} \right\|_{L^2}\\
	&\le C_{1}\left|\gamma-\varepsilon_0 \right|j^{-\zeta} \tau^{-\frac{3}{2}}(1+ \mathcal O (\tau^{-2} )), 
	\end{split}
	\end{equation}
	where $ C_{1} $ is a positive constant only depending on $\theta_0$ , $\delta $ and the volume  measure of $ \Cor_{r_{0}} $. Here $\delta$ and $ \Cor_{r_{0}} $ are given by \eqref{eq:innerdot} and \eqref{eq:coner0} respectively.
	 
In view of  \eqref{eq:vkr0} in Lemma~\ref{lem:v} and Lemma~\ref{lem:egv}, we know that 
	\begin{equation}\label{eq:333}
	\left| I_{6} \right| 
	\le Cj^{\beta}\tau^{-4}(1+ \mathcal O (\tau^{-2} )),
	\end{equation}
where $ \beta  $ is a positive constant,  $ C $ is a positive constant only depending on $\theta_0$ , $\delta $ and the volume  measure of $ \Cor_{r_{0}} $. 
	
	With the help of \eqref{eq:NonvanishIntegral} in Lemma~\ref{lem:Nonvanish}, we can obtain that 
	\begin{equation}\label{eq:334}
	\left|(\gamma-\varepsilon_0)\int_{\Cor_{r_{0}}} \mathbf{E}_{\mathbf  g_{j}}(\mathbf 0)\cdot \EM{V}\mathrm d \mathbf x \right| \ge \left|\gamma-\varepsilon_0 \right| \left|\mathbf{E}_{\mathbf g_{j}}(\mathbf 0)\cdot \mathbf p \right| C_{\mathcal {K}}(1+\frac{k^2}{\tau^2})^{-3/2}\tau^{-3}+\mathcal{O}(\tau^{-1}e^{-\frac{1}{2}r_{0}\tau \delta }),
	\end{equation}
		where  the positive number $C_{\mathcal K}$ is  independent  of $\tau$. 
    Due to $\gamma \neq  \varepsilon_0$,  choosing  $ \tau=j^{a} $ with $ a\in (\beta, \frac{2}{3}\zeta ) $, by virtue of  \eqref{eq:lem26}, \eqref{eq:331},  \eqref{eq:333} and \eqref{eq:334},  from \eqref{eq:equality} we derive that 
    \begin{equation}\label{eq:335}
    \begin{split}
    \left|\mathbf{E}_{\mathbf g_{j}}(\mathbf 0)\cdot \mathbf p \right| C_{\mathcal K}(1+\frac{k^2}{j^{2a} })^{-3/2}j^{-3a} &\le C j^{\beta}j^{-4a}(1+ \mathcal O (  j^{-2a} ))+C_{1}\left|\varepsilon-\varepsilon_0 \right|j^{-\zeta} j^{-\frac{3a}{2}}\\
    &\quad\times (1+ \mathcal O (  j^{-2a} ))  +{\mathcal O}( e^{-\delta r_0\tau } )
        \end{split}
    \end{equation}
    as  $\tau \rightarrow +\infty $. Multiplying $j^{3a}$ on both sides of \eqref{eq:335}, letting $j\rightarrow +\infty$, by noting $ a\in (\beta, \frac{2}{3}\zeta ) $,  we conclude that
    \begin{equation}\label{eq:246}
    \lim_{j\to\infty} \mathbf{E}_{\mathbf g_{j}}(\mathbf 0)\cdot \mathbf p=0.
    \end{equation}
   According to Lemma \ref{lem:vani yang} and \eqref{eq:246}, 
   we can deduce that
	\begin{equation}\nonumber
	\lim_{j\to\infty} \mathbf{E}_{\mathbf g_{j}}(\mathbf 0)=\mathbf 0. 
	\end{equation}
	Therefore, using \eqref{eq:eeg} we can prove \eqref{eq:thm1}.
	
	Finally, we can prove \eqref{eq:thm2} in a similar manner by choosing  $ \EM{V} $ and $ \EM{W} $ as \eqref{eq:vw2}.
	
	The proof is complete.
\end{proof}

	\section*{Acknowledgements}
The work of H. Diao is supported by a startup fund from Jilin University and NSFC No. 1211101002. The work of H. Liu is supported by the Hong Kong RGC General Research Funds (projects 12302919, 12301420 and 11300821) and the NSFC/RGC Joint Research Fund (project N\_CityU101/21).

\end{document}